\documentclass[11pt]{article}
\usepackage{amssymb,amsmath,mathtools}
\usepackage{color, xcolor}

\usepackage[OT2, T1]{fontenc} 
\usepackage[utf8]{inputenc} 

\usepackage[numbers,sort]{natbib}

\usepackage[mathlines]{lineno} 

\usepackage{centernot}

\usepackage[shortlabels]{enumitem}

\newcommand{\nequiv}{\not \equiv}

\usepackage[italian, english]{babel} 

\usepackage{csquotes} 

\usepackage
[draft=false,setpagesize=false,pdfstartview=FitH,
colorlinks=true,linkcolor=blue,citecolor=magenta,
pagebackref=true,
bookmarks=true, 
urlcolor=gray!65!black, 
]
{hyperref}

\usepackage[hyperpageref]{backref} 

\usepackage{footnotebackref} 

\numberwithin{equation}{section} 

\allowdisplaybreaks 

\def\claim#1.{\noindent {\bf #1.}}

\def\myflushright#1{{\unskip\nobreak\hfil\penalty50\hskip2em\hbox{}\nobreak\hfil%
#1\parfillskip=0pt\finalhyphendemerits=0\par}}

\def\bull{\vrule height 1.8ex width 1.0ex depth .1ex }
\def\QED{\ifmmode\eqno\hbox{$\bull$}\else\myflushright{\hbox{$\bull$}}\fi}

\newtheorem{Theorem}{Theorem}[section]
\newtheorem{Proposition}[Theorem]{Proposition}
\newtheorem{Corollary}[Theorem]{Corollary}
\newtheorem{Lemma}[Theorem]{Lemma}
\newtheorem{Remark}[Theorem]{Remark}
\newtheorem{Definition}[Theorem]{Definition}

\newcommand{\R}{\mathbb{R}}
\newcommand{\N}{\mathbb{N}}

\newcommand{\mc}[1]{\mathcal{#1}}

\newcommand{\wto}{\rightharpoonup}

\newcommand{\dx}{\mathrm{d}x}
\newcommand{\dy}{\mathrm{d}y}

\newcommand{\vertiii}[1]{{\left\vert\kern-0.25ex\left\vert\kern-0.25ex\left\vert #1 
		\right\vert\kern-0.25ex\right\vert\kern-0.25ex\right\vert}} 

\def\abs#1{|#1|}
\def\pabs#1{\left|{#1}\right|}
\def\norm#1{\|#1\|}
\def\intRN{\int_{\R^N}}

\def\eps{\varepsilon}

\newcommand{\sgn}{{\rm sgn}}

	\definecolor{azure(colorwheel)}{rgb}{0.0, 0.5, 1.0}
 	
\usepackage{scalerel}
\usepackage{tikz}
\usetikzlibrary{svg.path}
\definecolor{orcidlogocol}{HTML}{A6CE39}
\tikzset{
 orcidlogo/.pic={
 \fill[orcidlogocol] svg{M256,128c0,70.7-57.3,128-128,128C57.3,256,0,198.7,0,128C0,57.3,57.3,0,128,0C198.7,0,256,57.3,256,128z};
 \fill[white] svg{M86.3,186.2H70.9V79.1h15.4v48.4V186.2z}
 svg{M108.9,79.1h41.6c39.6,0,57,28.3,57,53.6c0,27.5-21.5,53.6-56.8,53.6h-41.8V79.1z M124.3,172.4h24.5c34.9,0,42.9-26.5,42.9-39.7c0-21.5-13.7-39.7-43.7-39.7h-23.7V172.4z}
 svg{M88.7,56.8c0,5.5-4.5,10.1-10.1,10.1c-5.6,0-10.1-4.6-10.1-10.1c0-5.6,4.5-10.1,10.1-10.1C84.2,46.7,88.7,51.3,88.7,56.8z};
 }
}
\newcommand\orcidicon[1]{\href{https://orcid.org/#1}{\mbox{\scalerel*{
\begin{tikzpicture}[yscale=-1,transform shape]
\pic{orcidlogo};
\end{tikzpicture}
}{|}}}}

\begin{document}

\title
{A Pohožaev minimization for normalized solutions: \\ fractional sublinear equations of logarithmic type}

\author{
		\\Marco Gallo 
		\orcidicon{0000-0002-3141-9598}
		\\ \normalsize{Universit\`{a} Cattolica del Sacro Cuore}
		\\ \normalsize{Dipartimento di Matematica e Fisica}
		\\ \normalsize{Via della Garzetta 48, 25133 Brescia, Italy}
		\\ \normalsize{\href{mailto:marco.gallo1@unicatt.it}{marco.gallo1@unicatt.it}}
		\\
		\\Jacopo Schino
		\orcidicon{0000-0003-2995-8637}
		\\ \normalsize{University of Warsaw}
		\\ \normalsize{Faculty of Mathematics, Informatics and Mechanics}
		\\ \normalsize{ul. Banacha 2, 02-097 Warsaw, Poland}
		\\ \normalsize{\href{mailto:j.schino2@uw.edu.pl}{j.schino2@uw.edu.pl}}
		\\
	}

\date{}

\maketitle


\begin{abstract}
In this paper, we search for normalized solutions to a fractional, nonlinear, and possibly strongly sublinear Schrödinger equation
$$(-\Delta)^s u + \mu u = g(u) \quad \hbox{in $\mathbb{R}^N$},$$
under the mass constraint $\int_{\mathbb{R}^N} u^2 \, \dx = m>0$; here, $N\geq 2$, $s \in (0,1)$, and $\mu$ is a Lagrange multiplier.
We study the case of $L^2$-subcritical nonlinearities $g$ of Berestycki--Lions type, without assuming that $g$ is superlinear at the origin, which allows us to include examples like a logarithmic term $g(u)= u\log(u^2)$ or sublinear powers $g(u)=u^q-u^r$, $0<r<1<q$. 
Due to the generality of $g$ and the fact that the energy functional might be not well-defined, we implement an approximation process in combination with a Lagrangian approach and a new Pohožaev minimization in the product space, finding a solution for large values of $m$. 
In the sublinear case, we are able to find a solution for each $m$.
Several insights on the concepts of minimality are studied as well. We highlight that some of the results are new even in the local setting $s=1$ or for $g$ superlinear.
\end{abstract}

\small


\medskip


%
%

\normalsize

{\hypersetup{linkcolor=black}
\tableofcontents
}

\section{Introduction}

\indent

In this paper, we study the fractional $L^2$-subcritical 
equation
\begin{equation}\label{eq_general_mu_fixed}
(-\Delta)^s u + \mu u = g(u) \quad \hbox{ in $\R^N$}
\end{equation}
under the mass constraint $\intRN u^2 \, \dx = m>0$; here, $s \in (0,1)$, $N>2s$, $\mu \in \R$ is a suitable Lagrange multiplier (part of the unknowns), and $g$ is a generic Berestycki--Lions \cite{BeLi83I} nonlinearity which is allowed to be (negatively) \emph{sublinear} at the origin (see \hyperref[g0]{\rm{(g0)}}--\hyperref[g4]{\rm{(g4)}} below). In particular, $G(t) := \int_0^t g(\tau) \, \mathrm{d}\tau$ is allowed to have a growth of the type
\begin{equation*}
 -|t|^{1+\theta} - |t|^{2^*_s} \lesssim
G(t) \lesssim |t|^2 + |t|^{\bar{p}+1}\quad \hbox{for $t \in \R$},
\end{equation*}
where $\theta \in [0,1)$, 
$2^*_s := \frac{2N}{N-2s}$ is the fractional Sobolev critical exponent, and
\begin{equation}\label{eq_p}
\bar{p}:= 1 + \frac{4s}{N};
\end{equation}
the number $\bar{p}+1$ is often referred to as the \textit{$L^2$-critical exponent}. 
This case includes, for example, the well-known \emph{logarithmic} nonlinearity $g(t) = t \log(t^2)$
(here $\theta \in [0,1)$ can be arbitrarily chosen), or \emph{sublinear-power}
-type nonlinearities 
$g(t) = |t|^{q-1}t - |t|^{r-1}t, \quad 0<r<q <\bar{p}$
(here $\theta =r$), 
but also strongly sublinear nonlinearities like $g(t) \sim \frac{\sgn(t) 
}{\log(|t|)} \quad \hbox{as $t\to 0$}$
(here $\theta=0$). 

Exploiting variational methods, a key fact of these nonlinearities is that the related energy functional $K\colon H^s(\R^N) \to \R \cup \{+\infty\}$, 
\begin{equation}\label{eq_def_energ_func_sfer}
K(u):=
\frac{1}{2}\intRN |(-\Delta)^{s/2}u|^2 \, \dx - \intRN G(u) \, \dx,
\end{equation}
is only lower semicontinuous, while it is well-posed in $ H^s(\R^N) \cap L^1(\R^N)$, which we highlight is not a reflexive space. Here $H^s(\R^N)$ denotes the usual fractional Sobolev space.


\subsection{Some literature}
\label{subsec_literature}


\indent


When $s=1$, 
the \emph{frequency} $\mu \in \R
$ 
is fixed, and $g(t)=t \log(t^2)$, the \emph{logarithmic Schrödinger equation} \cite{
BBM76, BBM79} 
\begin{equation}\label{eq_intr_log_s=1}
-\Delta u + \mu u = u \log(u^2) \quad \hbox{ in $\R^N$}
\end{equation}
has received much attention due to its great relevance in applied sciences, such as atomic physics, 
high-energy cosmic rays, Cherenkov-type shock waves, 
quantum hydrodynamical models, 
and many others; we refer to 
\cite{
	GLN10, Car22}
and the references therein. 
A key property of \eqref{eq_intr_log_s=1} is, indeed, the \emph{tensorization property} (or \emph{separability of noninteracting subsystems}): namely, if $u_i$ are solutions of \eqref{eq_intr_log_s=1} in $\R^{N_i}$ for $i=1,2$ with frequency $\frac{\mu}{2}$, then $
(x_1,x_2)\mapsto u_1(x_1)u_2(x_2)$ is a solution of \eqref{eq_intr_log_s=1} in $\R^{N_1} \times \R^{N_2}$.

In order to handle the singularity at the origin, several techniques have been developed: a first approach, applied by \cite{Caz83},
shows that the functional $K$ in \eqref{eq_def_energ_func_sfer} is of class $\mc{C}^1$ 
on the subspace 
$$X:= \left\{ u \in H^1(\R^N) \;\middle|\; \intRN u^2 \log(u^2) \, \dx >-\infty 
\right\};$$
see also \cite{
TaZh17,Shu19}, and \cite{BSV87} for general doubling nonlinearities.
We observe that $X \subsetneq H^1(\R^N)$, as shown by 
$u(x):= (|x|^{N/2} \log(|x|))^{-1}$ for $|x|$ large.
%
%
%
%
Finally, in \cite{LuNi19}, 
a suitable norm is introduced in order to force the positive solutions far from the singularity.
%
%



A second approach has been developed in \cite{DMS14, SqSz15}: 
here, the authors apply the \emph{nonsmooth analysis} theory of 
\cite{DeZa96, Szu86}
and work with
a suitable subdifferential of the action functional, possibly splitting it into a smooth part and a convex lower semicontinuous one. 

Different approaches rely instead on some \emph{approximation schemes}.
The idea that $\frac{t^{2+\delta}-t^2}{\delta} \to t^2 \log(t)$ as $\delta \to 0^+$ was formalized in the framework of PDEs by \cite{WaZh19} 
(see also \cite{GMS26}, 
\cite[Section 7]{Car22}), 
and then this idea was applied in \cite{ZhWa20} 
to get existence of solutions. 
See also \cite{GLN10}, 
where the regularization is of the form $\log(t+\delta)$, 
and \cite{AlJi20}, 
where the authors 
approximate the domain of the problem with expanding balls. 




We mention that the problem of the singularity 
of the logarithm arises also in the study of planar Choquard equations 
\cite{CiWe16,LRTZ22}.

Let us move, now, to the case of sublinear powers. Although nonlinear (fractional) Schrödinger equations with combined powers have recently gained the mathematical community's attention \cite{DiFe19,TVZ07, LuZh20, Soa20s} (see also \cite{BdS} for a more general nonlinear term that behaves similarly), such work usually concerns superlinear exponents. Instead, as
regards sublinear-power nonlinearities (again, $s=1$ and $\mu$ fixed),
most of the work focuses on problems of the type
$$-\Delta u + V(x) u = a(x) |u|^{q-1}u + h(x) \quad \hbox{ in $\R^N$},$$
$q \in (0,1)$, where suitable assumptions on $V, a, h$ 
are considered in order to ensure the existence of solutions (notice, in fact, that when $V\equiv a \equiv 1$ and $h \equiv 0$, no solution exists, 
see e.g. 
\cite[Lemma 2.3]{BeOu09}). 
When $a \not\in L^{\frac{2}{1-q
}}(\R^N)$, 
the analysis of the problem is usually set in the Banach space $H^1(\R^N) \cap L^{q+1}(\R^N)$;
see 
\cite{Teh07, BeOu09, 
BOR15}. 


%
The only paper we know that deals with \emph{general nonlinearities} in the spirit of \cite{BeLi83I} is \cite{Med21}: 
here, the author studies
$$-\Delta u 
= f(u) \quad \hbox{ in $\R^N$},$$
where the assumptions about $f$ include the case
$\lim_{t \to 0} \frac{f(t)}{|t|^{2^*-2}t}=-\infty$. The author employs a 
perturbation argument at the origin
and obtains a ground state for the equation. 


When $s=1$ and $\mu$ is free, very few results appear in the literature about the existence of normalized solutions $\intRN u^2 \, \dx = m > 0$. When $g(t)=t \log(t^2)$, we see that, considered a nontrivial solution $u_{1}$ of \eqref{eq_intr_log_s=1} for $\mu=1$, then 
defining $u:= \alpha_{m} u_{1}$ with $\alpha_{m}^2 := m (\intRN u_1^2 \, \dx)^{-1}$, we have that $u$ satisfies $\intRN u^2 \, \dx = m$ and
\begin{equation}\label{eq_intr_scaling}
-\Delta u + \left( 1 + \log(\alpha_{m}^2)\right) u = u \log(u^2) \quad \hbox{in $\R^N$,}
\end{equation}
thus the existence of a normalized solution is always ensured; this was already noticed in \cite{BBM79}. 
In particular, 
the frequency 
$1 + \log(\alpha_{m}^{2})$ is negative when $m$ is sufficiently small and positive when $m$ is sufficiently large. 
Additionally, an explicit $L^2$-minimum can be found. 
Indeed, by the logarithmic Sobolev--Gross--Nelson inequality \cite{Gro75} there exists $C(N) \in \R$ such that 
%
\begin{equation}\label{eq_intr_optineq}
K(u) \geq - \frac12 \left( C(N) 
+ \log(m)\right) m
\end{equation}
for every $u \in H^1(\R^N)$ with $\intRN u^2 \, \dx = m$; moreover, we know that the equality is attained by (and only by, up to a translation) the \emph{Gausson} $u(x)=\gamma_m e^{-\frac{1}{2} |x|^2}$ with $\gamma_m:= \sqrt{m} 
\left(\intRN e^{-|x|^2} \, \dx\right)^{-1/2}
$. 
This procedure can be found, for example, in \cite{BBM76, Caz83}. 
We see that the sign of the minimal energy over the $L^2$-sphere depends on the size of $m$ (and this is connected with the sign of the Lagrange multiplier, see Proposition \ref{prop_intr_L2min_PN}).


When $g$ does not enjoy scaling properties of this type, we mention 
the 
recent papers 
\cite{ShYa23, AlJi24, ZhZh24, MeSc24} 
(see also \cite{DiJiPu25, JiYa26} for results in the quasilinear setting). 
In \cite{ShYa23} and \cite{AlJi24}, 
the authors study, respectively, logarithm-plus-power and nonautonomous problems by restricting to 
the subspace $X$.
%
%
%
In \cite{ZhZh24}, via a direct minimization, the authors study a similar equation, but with a more general $g$ satisfying a monotonicity assumption: 
they find an $L^2$-minimum for every mass $m$, which, in particular, has positive energy for $m$ small 
and negative energy for $m$ large. 
Finally, by exploiting the perturbation techniques developed in \cite{Med21}, the authors in \cite{MeSc24} succeed in removing the monotonicity condition on $g$, and find an $L^2$-minimum for 
$m$ 
large, with a positive multiplier $\mu>0$. 


When $s \in (0,1)$, the 
%
fractional logarithmic equation
\begin{equation}\label{eq_intr_frac_log}
(-\Delta)^s u + \mu u = u \log(u^2) \quad \hbox{in $\R^N$}
\end{equation}
appears in the study of nonlinear quantized boson fields 
\cite{Ros69}. 
Mathematically, several papers succeeded in generalizing the results found in the local case for sublinear nonlinearities, exploiting and adapting to the nonlocal case the aforementioned techniques, see e.g.
\cite{DSZ14, Ard17, JiXu22} 
for the logarithmic equation, 
and \cite{FaZh23}
for 
logarithm-plus-power nonlinearities.
In particular, we mention \cite{Iko20}: 
here, the author, with a perturbation argument similar to the one in \cite{Med21}, 
finds infinitely many solutions with a general nonlinear term. 


Anyway, all the previous results deal with the unconstrained problem with fixed frequency.
When $g(t)=t \log(t^2)$, the search for a normalized solution can be pursued as above by scaling (see \eqref{eq_intr_scaling}); on the other hand, the search for an $L^2$-minimum through an optimal logarithmic inequality \eqref{eq_intr_optineq} has issues.
Indeed, although such an inequality still holds in the fractional framework 
\cite{CoTa05, ChRu24}, 
i.e., there exists $C(N,s) \in \R$ such that 
$$K(u) \geq - \frac12 \left( C(N,s) 
+\log(m)\right) m$$
for every $u \in H^s(\R^N)$ with $\intRN u^2 \, \dx = m$,
it is not known if the equality is attained or not.
Normalized solutions to fractional equations with logarithmic nonlinearities have only recently begun to be investigated; see, e.g.,  \cite{LiZh25,LvLi24}, where semiclassical problems are studied.


The goal of this paper is to find normalized solutions with minimality properties to fractional sublinear equations under general assumptions about the nonlinear term.
When treating such general terms $g$, as previously commented, the splitting of the functional does not seem an option, as well as the definition of a space similar to $X$.
To obtain the existence of an $L^2$-minimum we rely, thus, on the perturbation technique employed by 
\cite{Med21, MeSc24}: 
this perturbation takes the form $g_{\eps}(t):= g_+(t) - \varphi_{\eps}(t) g_-(t)$, where
$\varphi_{\eps}(t) = |t|/\eps$ for $|t|\leq \eps$ (see Section \ref{s-pert} for details). Additionally, we obtain several properties that are new even for $s=1$ and $g$ superlinear. 

\subsection{Existence results}

\indent

In the present paper, we study 
the following problem
\begin{equation}\label{e-main}
\begin{cases}
\displaystyle (-\Delta)^s u + \mu u = g(u) \quad \hbox{in $\R^N$},\\
\displaystyle \intRN u^2 \, \dx = m,\\
(\mu,u) \in \R \times H^s(\R^N),
\end{cases}
\end{equation}
with $s \in (0,1)$, $N \ge 2$, $m > 0$, and $g \colon \R \to \R$ satisfying (set $G(t) := \int_0^t g(\tau) \, \mathrm{d}\tau$)
\begin{itemize}
	\item[(g0)] \label{g0}
	$g$ is continuous and $g(0) = 0$,
	\item[(g1)] \label{g1}
	$\displaystyle \limsup_{t \to 0} \frac{g(t)}{t} \leq 0$,
	\item[(g2)] \label{g2}
	$\displaystyle \limsup_{|t| \to +\infty} \frac{|g(t)|}{|t|^{2_s^*-1}} < +\infty$, 
	\item[(g3)] \label{g3}
	$\displaystyle \limsup_{|t| \to +\infty} \frac{g(t)}{|t|^{\bar{p}-1}t} \leq 0$,
	\item[(g4)] \label{g4}
	there exists $t_0 \ne 0$ such that $G(t_0) > 0$.
\end{itemize}
%
Here and in what follows, we denote $g^{\pm}:=\max\{\pm g,0\}$ the standard positive and negative parts of $g$, and
\begin{align}
&g_+(t)t:=(g(t)t)^+, \quad G_+(t):= \int_0^t g_+(\tau) \, \mathrm{d}\tau, \label{eq_intr_g+}
\\&g_-(t)t:=(g(t)t)^-, \quad G_-(t):= \int_0^t g_-(\tau) \, \mathrm{d}\tau; \label{eq_intr_g-}
\end{align}
see Section \ref{s_prelim}. 
With these notations, we may rewrite \hyperref[g1]{\rm{(g1)}} and \hyperref[g3]{\rm{(g3)}} as
\begin{itemize}
	\item[(g1)] $\displaystyle \lim_{t \to 0} \frac{g_+(t)}{t} = 0$, 
	\item[(g3)] $\displaystyle \lim_{|t| \to +\infty} \frac{g_+(t)}{|t|^{\bar{p}-1}t} = 0$.
\end{itemize}
We further notice that, when \hyperref[g3]{\rm{(g3)}} holds, then \hyperref[g2]{\rm{(g2)}} actually means 	$\limsup_{|t| \to +\infty} 
\frac{g_-(t)}{|t|^{2_s^*-2}t} < \infty$.

\medskip


Before presenting our results, let us recall that in \cite[Proposition 4.1]{ChWa13} 
and \cite[Proposition 1.1]{BKS17} 
(see also \cite[Corollary 6.4]{CGT24} for a different approach), 
when $g$ satisfies Berestycki--Lions assumptions, the following \emph{Poho\v{z}aev identity} holds
\begin{equation}\label{eq_pohoz_frac}
\frac{1}{2^*_s} \int_{\R^N}|(-\Delta)^{s/2} u|^2 \, \dx + \frac{\mu}{2} \int_{\R^N} u^2 \, \dx - \intRN G(u) \, \dx = 0
\end{equation}
whenever 
$u$ is a solution of \eqref{e-main} and either $s \in (\frac{1}{2}, 1)$ or $g \in \mc{C}^{0,\alpha}_{loc}(\R^N)$ with $\alpha \in (1-2s,1)$: this result can be generalized to our assumptions whenever $\intRN G(u) \, \dx$ is assumed a priori finite; however, we see that the condition on the Hölder regularity of $g$ forces a restriction on the singularity of $g_-$ at the origin.
It is not known if this identity holds for general continuous nonlinearities $g$ and general values of $s \in (0,1)$. 
Additionally, we observe that, due to the possible sublinearity of $g$, for a general solution $u$ of equation \eqref{e-main} we cannot state a priori that 
it satisfies the \emph{Nehari identity} 
\begin{equation}\label{eq_nehar_id}
 \intRN |(-\Delta)^{s/2} u|^2 \, \dx + \mu \intRN u^2 \, \dx - \intRN g(u) u \, \dx =0.
 \end{equation}

To state the main theorem, we make use of a 
\emph{Lagrangian formulation} of the problem (in the spirit of \cite{HiTa19, CGT21N}) and introduce, in addition to \eqref{eq_def_energ_func_sfer}, a Lagrangian functional over the product space: $I^m \colon \R \times H^s(\R^N) \to \R \cup \{+\infty\}$,
\begin{equation}\label{eq_def_Im}
I^m(\mu,u) := \frac12 \int_{\R^N}|(-\Delta)^{s/2} u|^2 \, \dx + \frac{\mu}{2} \left(\int_{\R^N} u^2 \, \dx - m\right) - \intRN G(u) \, \dx.
\end{equation}
Notice that critical points of $I^m$ are (formally) 
solutions to \eqref{e-main} and $I^m(\mu,u) = K(u)$ if $\int_{\R^N} u^2 \, \dx = m$.

From now on, to state our results we will use the following convention.

\smallskip

\textbf{Convention:} whenever we deal with a pair $(\mu,u)$ (or $(\lambda,u)$ with $e^\lambda = \mu$, see \eqref{eq_identif_lambda_mu}) and we say that such a pair is radially symmetric, nonnegative, does not change sign, etc., \emph{we are referring to $u$}.

\smallskip

We introduce the following notations, standing for the $L^2$-sphere in $ H^s(\R^N)$ and the Poho\v{z}aev set in the product space $\R \times H^s(\R^N)$:
\begin{equation}\label{eq_def_Sm}
\mc{S}_m:=\big\{u \in H^s(\R^N) \mid |u|_2^2=m\big\},
\end{equation}
where $|\cdot|_q$ is the usual Lebesgue norm, $1 \le q \le \infty$, and 
\begin{equation}\label{eq_def_mcP_pos}
	\mc{P} := \left\{(\mu, u) \in \R \times H^s(\R^N) \mid u \ne 0 \text{ and } P(\mu,u) = 0\right\},
\end{equation}
where $P \colon \R \times H^s(\R^N) \to \R \cup \{+\infty\}$ is the Poho\v{z}aev functional%
\footnote{We highlight that the terminology ``Poho\v{z}aev functional'', in the framework of normalized solutions, is often referred to a suitable combination of the Poho\v{z}aev and Nehari functionals, which allows getting rid of the Lagrange multiplier $\mu$; such a functional appears, in particular, in the framework of $L^2$-supercritical problems. Here, instead, we focus on the proper Poho\v{z}aev identity, considering a two-variable functional.}
\begin{equation}\label{eq_def_Poh_pos}
	P(\mu,u) := \intRN |(-\Delta)^{s/2}u|^2 \, \dx + 2_s^* \intRN \frac{\mu
	}{2} u^2 - G(u) \, \dx.
\end{equation}
It is clear that $(\mu,u) \in \mc{P}$ if and only if it satisfies the Poho\v{z}aev identity \eqref{eq_pohoz_frac}.
Moreover, $\mc{P} \ne \emptyset$ because, fixed $u \in H^s_\textup{rad}(\R^N) \setminus \{0\}$ with $\int_{\R^N} G(u) >-\infty$, $P(\cdot,u)$ is a continuous function such that $\lim_{\mu \to \pm \infty} P(\mu,u) = \pm \infty$.
Next, we set 
$$\mc{S}_m^{\textup{rad}}:=\mc{S}_m \cap H^s_{\textup{rad}}(\R^N), \quad \mc{P}^{\textup{rad}}:=\mc{P} \cap \bigl(\R \times H^s_{\textup{rad}}(\R^N)\bigr),$$
where $H^s(\R^N)$ is the classical fractional Sobolev space, cf. \eqref{def_Ds_Hs}, and $H^s_{\textup{rad}}(\R^N)$ is its subspace of radially symmetric functions,
and
\begin{equation}\label{eq_def_km_em}
\kappa^m 
:= \inf_{\mc{S}_m^{\textup{rad}}} K, \qquad d^m:= \inf_{\mc{P}^{\textup{rad}}} I^m;
\end{equation}
here $K$ and $I^m$ are defined in \eqref{eq_def_energ_func_sfer} and \eqref{eq_def_Im}.
We further define 
\begin{equation}\label{eq_def_mubar0}
\overline{\mu}_0:= \sup_{t\ne0} \frac{G(t)}{t^2/2} 
;
\end{equation}
notice that if \hyperref[g4]{\rm{(g4)}} holds, then $
\overline{\mu}_0 \in (0,+\infty]$.
%
%
%
%
%
%
%
Moreover, we introduce the $L^2$-ball and the corresponding infimum:
\begin{equation}\label{eq_def_Dm_ellm}
\mc{D}_m:= \left\{ u \in H^s(\R^N) \mid |u|_2^2 \leq m\right\}, \quad \mc{D}_m^{\textup{rad}}:=\mc{D}_m \cap H^s_{\textup{rad}}(\R^N), \quad \ell^m:= \inf_{\mc{D}_m^{\textup{rad}}} K;
\end{equation}
we mention that the use of the $L^2$-ball in normalized problems goes back to 
\cite{BiMe21}, while it was utilized for the first time in the $L^2$-subcritical regime in \cite{Sch22} and in a mixed regime in \cite{BdS}.
We can thus state our 
first 
theorem, see Definition \ref{def_weak_sol} for the notion of solution. 
\begin{Theorem}[Existence for large masses]\label{t-main}
If \hyperref[g0]{\rm{(g0)}}--\hyperref[g4]{\rm{(g4)}} hold, then there exists $m_0 \geq 0$ such that for every $m > m_0$ there exists $(\mu_0,u_0) \in (0,\infty) \times \mc{S}_m^\textup{rad} 
$ 
with the following properties:
\begin{enumerate}[(i)]
	\item $(\mu_0,u_0)$ is a 
	solution to \eqref{e-main}, and it
	%
		satisfies the Nehari identity \eqref{eq_nehar_id};
	
	\item $u_0$ is a minimum 
	(with Lagrange multiplier $\mu_0$) on the $L^2$-sphere and the $L^2$-ball, i.e., $K(u_0) = \kappa^m = \ell^m 
	<0$;
	moreover, we have the formula
	\begin{equation}\label{eq_legendre_intr}
 \kappa^m 
 = \inf_{-\infty< \mu < \overline{\mu}_0} \left( a(\mu) - \mu \frac{m}{2}\right) = a(\mu_0) - \mu_0 \frac{m}{2},
 \end{equation}
 where 
$a(\mu)$ is the least action related to the problem with fixed frequency $\mu$, see 
\eqref{eq_def_amu}. 
	
	\item $(\mu_0, u_0)$ is a Poho\v{z}aev minimum on the product space, i.e., $(\mu_0,u_0) \in \mc{P}^{\textup{rad}} 
	$ and $I^m(\mu_0,u_0) = d^m$; 
in particular, $(\mu_0, u_0)$ satisfies the Poho\v{z}aev identity \eqref{eq_pohoz_frac}.
%


\end{enumerate}
%
As a consequence of the above relations, we have $\kappa^m 
=d^m$.
\end{Theorem}

Under suitable additional assumptions about $g$, we highlight that the found solution can be chosen nonnegative and radially nonincreasing, see Proposition \ref{p-sym_unpert}.


\begin{Remark}
We observe the following facts.
\begin{itemize}
	\item The mass threshold $m_0$ is given 
	by formula \eqref{eq_meps_m0} and it ensures that $\kappa^m < 0$ if $m > m_0$. A detailed analysis on mass thresholds for fractional Schr\"odinger equations is provided in \cite{BieSch}. 
	\item Notice that, thanks to our method, we are able to state that the found Lagrange multiplier $\mu_0$ is positive and the energy $\kappa^m$ is negative. 
	At the same time, the found solution $(\mu_0, u_0)$ minimizes $I^m$ in the whole Poho\v{z}aev space over $\R\times H^s_{\textup{rad}}(\R^N)$, not only over $(0,+\infty)\times H^s_{\textup{rad}}(\R^N)$. 
	\item In Theorem \ref{t-main}, the solution is found for $m$ sufficiently large. In Theorem \ref{t-small-mass} below, we find a solution for each $m>0$; on the other hand, in that case, we have no information on the sign of the Lagrange multiplier or that of the energy, 
	nor do we know whether 
	the solution minimizes the energy over $\mc{D}_m^\textup{rad}$ (which is coherent with Proposition \ref{prop_uguag_not_sym} below).
	\item We highlight that we cannot state that $u_0$ is a critical point of $I^m(\mu_0, \cdot)$ due to the lack of regularity of the functional, but only that $I^m(\mu_0, \cdot)$ is differentiable at $u_0$ along test functions in $H^s(\R^N) \cap L^1(\R^N)$, cf. Lemma \ref{lem_diff_along_L1}.
In particular, 
 the validity of the Nehari identity in Theorem \ref{t-main} is not straightforward; 
in this case, thanks to our method, we are able to obtain such an identity.
\item We observe that \eqref{eq_legendre_intr} is a ``Legendre transform''-type formula which relates the $L^2$-minimum value to the Poho\v{z}aev minimum value (with fixed frequency); similar formulas were already obtained in \cite{DST23, HiTa19, Gal23}.
\item We notice that, fixed $\rho \in (0, \overline{\mu}_0)$, we can find a threshold $m_{\rho}\geq m_0$ such that, for each $m>m_{\rho}$ there exists $(\mu_{\rho}, u_{\rho} ) \in ( \rho, +\infty) \times \mc{S}_m^\textup{rad}$ 
solution of \eqref{e-main}. Indeed, it is sufficient to apply Theorem \ref{t-main} to the nonlinearity $g_{\rho}(t):= g(t) - \rho t$, which satisfies \hyperref[g0]{\rm{(g0)}}--\hyperref[g4]{\rm{(g4)}}. 
\end{itemize}
\end{Remark}


\smallskip

Unlike \cite{MeSc24}, 
in our fractional setting, two main difficulties come into play. 
The first is due to the nonlocality of the problem: 
the fact that the operator $(-\Delta)^s$ does not preserve the supports of functions
creates indeed problems in a 
direct minimization 
as in \cite{MeSc24}. 
To deal with this issue, we handle the perturbed problem through a Lagrangian formulation \eqref{eq_def_Im} as done in \cite{HiTa19, CGT21N}: 
this approach has additional advantages, like showing some geometrical properties of the solutions, such as the minimality of $\mc{P}^{\textup{rad}}$, formula \eqref{eq_legendre_intr}, and the form of the threshold $m_0$; 
this last, in turn, allows us to show the existence of a solution for every mass $m$ (see Theorem \ref{t-small-mass} below).

The second difficulty is due to the absence of a Poho\v{z}aev identity: as we already pointed out, it 
is an open problem to determine if the Poho\v{z}aev identity holds for general continuous nonlinearities $g$ and general values of $s \in (0,1)$. Under our assumptions \hyperref[g0]{\rm{(g0)}}--\hyperref[g4]{\rm{(g4)}}, thus, we are not able to use this identity. 
As already mentioned, we first study the superlinear problem (with an $\eps$-perturbation) and pass to the limit: 
by construction, the solutions $u_{\eps}$ that we build satisfy the Poho\v{z}aev identity, which gives some information to start with. 
On the other hand, in the limiting process, we get $u_{\eps} \wto 
u_0$, where $
u_0$ can be proved to be a solution, but with no a priori minimality properties: 
in \cite{MeSc24}, it is thus crucially used, in order to pass to the strong limit and get the existence of an $L^2$-minimum, that 
the Poho\v{z}aev identity holds for arbitrary solutions.
This obstruction is here overcome by a finer analysis on the information of $u_{\eps}$, in particular, regarding the interplay between the minimization over the $L^2$-sphere and the $L^2$-ball.

To deal with the superlinear perturbed prescribed-mass problem, in the spirit of \cite{HiTa19, CGT21N}, 
we first study the free problem (i.e., when $\mu$ is fixed), following 
essentially 
\cite{ChWa13, BKS17}: 
here 
we revise and refine both the existence results contained in \cite{BKS17, CGT21N}.
Additionally, since we are interested in the existence of a single solution, we simplify the argument of \cite{CGT21N} by finding the solution through a \emph{direct minimization on the Poho\v{z}aev set} \eqref{eq_def_mcP_pos}: as a matter of fact, the minimality of the solution in \cite{CGT21N} was obtained as a by-product of a minimax approach (see Remark \ref{rem_multipl}), which requires a nontrivial deformation theory. Here we directly study 
the convergence of a (precisely chosen) minimizing sequence on $\mc{P}^{\textup{rad}}$: this procedure is totally unrelated to 
the perturbation framework (see Remark \ref{t-g*}), and it seems a new approach in literature in regards of product spaces, thus it has an interest of its own.

\smallskip

In Theorem \ref{t-main}, we show the existence of a minimizer for $m>0$ sufficiently large; when the mass is small, by looking at the pure logarithmic case we expect an $L^2$-minimum, possibly with a Lagrange multiplier 
and an energy with signs different from the ones in Theorem \ref{t-main}. This is the content of the next 
theorem 
(see also the more general Theorem \ref{thm_small_mass}). 
This result gives a fractional counterpart of \cite{ZhZh24}, as well as a generalization for $s=1$ by dropping the monotonicity therein.

\begin{Theorem}[Existence for all masses]\label{t-small-mass}
If \hyperref[g0]{\rm{(g0)}}--\hyperref[g4]{\rm{(g4)}} hold and, moreover,
\begin{equation}\label{eq_cond_meno_inf}
\lim_{t \to 0} \frac{g(t)}{t} = -\infty,
\end{equation}
then for every $m>0$ there exists a 
solution $(\mu_0,u_0) \in \R \times \mc{S}_m^\textup{rad}$ 
of \eqref{e-main} satisfying the Nehari and Poho\v{z}aev identities and such that
$$I^m(\mu_0, u_0) = K(u_0) = a(\mu_0) - \frac{\mu_0}{2} m = \kappa^m 
= d^m = \inf_{-\infty< \mu < \overline{\mu}_0} \left(a(\mu)-\frac{\mu}{2}m\right) .$$
	%
	%
	Finally, if $\kappa^m \leq 0$, 
	then $\mu_0>0$.
\end{Theorem}
We remark that in Theorem \ref{t-small-mass}, 
we are not able to determine the sign of the Lagrange multiplier $\mu_0$ or of the energy 
$\kappa^m$. 
To obtain Theorem \ref{t-small-mass}, we make use of a shifted problem, which still falls into the assumptions of Theorem \ref{t-main} (here the generality of \hyperref[g1]{\rm{(g1)}} is essential), and rely on the representation formula \eqref{eq_meps_m0} for the mass threshold $m_0$. In particular, in obtaining solutions for every $m>0$, \eqref{eq_cond_meno_inf} plays a totally different role than the condition $\lim_{t \to 0} G(t) / |t|^{\bar{p}+1} = +\infty$ (cf. Remark \ref{rem_multipl}), which is used to ensure that $\kappa^m = d^m$ is negative for every $m>0$; see \cite{CGT21N,JeLu22O,Sch22} and the references therein.

\begin{Remark}[On the multiplicity]
\label{rem_multiplicity}
We highlight that the present approach seems suitable for generalizations in the 
$L^2$-critical (see \cite{CGIT24}) and supercritical cases.
Moreover, by following the minimax ideas in \cite{CGT21N}, 
one may have the additional advantage to possibly obtain 
a multiplicity result: indeed, for each $\eps>0$, with the same techniques (and genus arguments, see Remark \ref{rem_multipl}) 
we can get the existence of multiple solutions $(u^{\eps}_k)_{k\in \N}$; on the other hand, at the moment, we are not able to exclude the possibility that these solutions collapse 
to the same solution of the original problem when $\eps \to 0^+$; we expect that some finer topological information on the genus of these solutions is needed.
\end{Remark}

\subsection{Qualitative results and examples}
\label{subsec_furth_results}


\indent

We present now some outcomes about the interplay among the various concepts of minimality we introduced, as well as some symmetry results. Most of them will not require \hyperref[g3]{\rm{(g3)}}. 
%
%
%

As highlighted in Subsection \ref{subsec_literature}, 
we recall that, in some cases, the Lagrange multiplier of the equation might be negative. In the following, thus, for $\rho \in [-\infty,+\infty)$ we define
\begin{equation} \label{eq_def_pvarpi}
\mc{P}_{(\rho)} := \left\{(\mu, u) \in (\rho,+\infty) \times H^s(\R^N) \mid u \ne 0 \text{ and } P(\mu,u) = 0\right\}
\end{equation}
and $\mc{P}_{(\rho)}^{\textup{rad}} := \mc{P}_{(\rho)} \cap \big( \R \times H^s_\textup{rad}(\R^N) \big)$. Observe that $\mc{P}_{(-\infty)} = \mc{P}$. 

In the following propositions, we are not claiming that the mentioned infima are attained or finite 
merely under \hyperref[g0]{\rm{(g0)}}--\hyperref[g2]{\rm{(g2)}}. On the other hand, 
minimizers can exist under different assumptions than those in Theorems \ref{t-main} or \ref{t-small-mass}: for instance, with proper hypotheses on $m$ and $g$, they exist when $\limsup_{|t| \to +\infty} 
g(t)/ (|t|^{\bar{p}-1}t) < +\infty$ (which is weaker than \hyperref[g3]{\rm{(g3)}}); see, e.g., \cite{LuZh20}.

Under suitable assumptions, the several notions of minimality we introduced coincide; notice that when $m>m_0$, we have $\kappa^m<0$, which is a particular case of such assumptions. 
\begin{Proposition}[Connections among infima]
\label{prop_uguag_not_sym}
Assume \hyperref[g0]{\rm{(g0)}}--\hyperref[g2]{\rm{(g2)}}.
We have 
\begin{equation}\label{eq_ineq_not_sym}
\inf_{\mc{S}_m} K = \inf_{\mc{P}} I^m.
\end{equation}
If $\rho 
\in \R$ and $\inf_{\mc{S}_m} K < -\frac{\rho}{2}m$, then
$$
\inf_{\mc{S}_m} K = \inf_{\mc{P}} I^m = \inf_{\mc{P}_{(\rho)}} I^m.$$
If $\rho \geq 0$ and $\inf_{\mc{D}_m} K < -\frac{\rho}{2}m$, then
$$\inf_{\mc{D}_m} K = \inf_{\mc{S}_m} K = \inf_{\mc{P}} I^m = \inf_{\mc{P}_{(\rho)}} I^m.$$
In particular, if $\inf_{\mc{D}_m} K <0$ we have 
$\inf_{\mc{D}_m} K = \inf_{\mc{S}_m} K.$

Finally, for every $\rho \in [-\infty,+\infty)$, we have
$$\inf_{\mc{P}_{(\rho)}} I^m = \inf\Big\{ \frac{s}{N} |(-\Delta)^{s/2} u|_2^2 -\frac{\mu}{2} m \;\Big|\; (\mu, u) \in (\rho, +\infty) \times H^s(\R^N), \, u \neq 0, \, P(\mu, u) \leq 0\Big\},$$
and $\inf_{\mc{P}_{(\rho)}} I^m < \frac{s}{N} |(-\Delta)^{s/2} u|_2^2 -\frac{\mu}{2} m$ if $P(\mu, u) < 0$.
\end{Proposition}
As shall be evident from the proof, the conditions about $m$, $\rho$, and the infimum of $K$ over $\mc{S}_m$ or $\mc{D}_m$ are needed to ensure that there exists $(\mu,u) \in \mc{P}_{(\rho)}$ such that $u \in \mc{D}_m$. 

We now show a generalization of \eqref{eq_legendre_intr}. For $\mu \in \R$, we denote by $J_\mu$ and $\mc{P}_\mu$ the action functional and the Poho\v{z}aev manifold related to the free-mass problem, see \eqref{eq_def_Jmu} and \eqref{eq_def_poh_set} respectively.

\begin{Proposition}[Legendre transform formula]
	\label{prop_intr_legendre}
	Assume \hyperref[g0]{\rm{(g0)}}--\hyperref[g2]{\rm{(g2)}} and let $\rho \in [-\infty,\overline{\mu}_0)$.
	Then,
	$$ \inf_{\mc{P}_{(\rho)}} I^m 
	= \inf_{\rho < \mu < \overline{\mu}_0} \left(\inf_{\mc{P}_{\mu}} J_{\mu}-\frac{\mu}{2}m\right).$$
	%
	%
\end{Proposition}


We move on to show that Poho\v{z}aev minima are indeed always weak solutions.

\begin{Proposition}[Poho\v{z}aev minima are solutions]\label{pr-Poho_sol}
Assume \hyperref[g0]{\rm{(g0)}}--\hyperref[g2]{\rm{(g2)}} and let $\rho \in [-\infty,+\infty)$. If $(\mu,u) \in \mc{P}_{(\rho)}$ satisfies $I^m(\mu,u) = \min_{\mc{P}_{(\rho)}} I^m$, then it is a solution to \eqref{e-main}.
 Moreover, if 
$\min_{\mc{P}_{(\rho)}} I^m < -\frac{\rho}{2}m$, then $I^m(\mu,u) = K(u)= \inf_{\mc{S}_m} K$.
\end{Proposition}

We show now that if $u$ is an $L^2$-minimum, then $u$ satisfies the Poho\v{z}aev identity; this is obtained without assuming extra regularity.
See also Proposition \ref{p-Nehari} for some results on the Nehari identity.



\begin{Proposition}[Poho\v{z}aev identity for $L^2$-minima] 
\label{prop_intr_L2min_PN}
Assume \hyperref[g0]{\rm{(g0)}}--\hyperref[g2]{\rm{(g2)}}.
Let $u \in \mc{S}_m$ be such that $K(u)= \inf_{\mc{S}_m} K$. 
Then, $u$ is a 
solution of \eqref{eq_general_mu_fixed} for some Lagrange multiplier $\mu \in \R$, and $(\mu,u)$ satisfies the Poho\v{z}aev identity.
Moreover, for every $\rho \in [-\infty, - \frac{2}{m} \inf_{\mc{S}_m} K)$, we have
$K(u) = I^m(\mu,u) = \inf_{\mc{P}_{(\rho)}} I^m$.
%
\end{Proposition}

The conclusion of Proposition \ref{prop_intr_L2min_PN} holds for $u \in \mc{D}_m \setminus \mc{S}_m$ such that $K(u) = \inf_{ \mc{D}_m } K$; in such a case, $\mu=0$. See Proposition \ref{pr:PohoMin}.
As a consequence of Propositions \ref{pr-Poho_sol} and \ref{prop_intr_L2min_PN} we have the following statement.

\begin{Corollary}\label{corol_equival_min_min}
Under \hyperref[g0]{\rm{(g0)}}--\hyperref[g2]{\rm{(g2)}}, a function $u\in H^s(\R^N)$ is an $L^2$-minimum with Lagrange multiplier $\mu \in \R$ if and only if $(\mu, u)$ is a Poho\v{z}aev minimum over $\mc{P}$.
\end{Corollary}


Finally, we present a fractional and possibly strongly sublinear counterpart of \cite{JeLu22O}: notice that here we deal with Poho\v{z}aev minima instead of classical ground state solutions, which are indeed equivalent notions in the case treated in \cite{JeLu22O}.

\begin{Proposition}[Least energy versus least action]
\label{prop_JJLU}
Assume \hyperref[g0]{\rm{(g0)}}--\hyperref[g2]{\rm{(g2)}}. 
\begin{itemize}
\item[(i)] Let $u \in \mc{S}_m$ be such that $K(u)= \inf_{\mc{S}_m} K$ with Lagrange multiplier $\mu \in \R$. Then $u \in \mc{P}_{\mu}$ is such that $J_{\mu}(u)=\inf_{\mc{P}_{\mu}} J_{\mu}$ 
and
\begin{equation}\label{eq_relat_kep_a}
\inf_{\mc{S}_m} K = \inf_{\mc{P}_{\mu}} J_{\mu} - \frac{\mu}{2} m.
\end{equation}
\item[(ii)]
Let $$\mu \in \Big\{ \nu \in \R \;\Big|\; \hbox{$\exists$ $v \in \mc{S}_m$ s.t. $K(v)= \inf_{\mc{S}_m} K$ with Lagrange multiplier $\nu$}\Big\}.$$ 
Let now $u \in \mc{P}_{\mu}$ be such that $J_{\mu}(u)=\inf_{\mc{P}_{\mu}} J_{\mu}$. Then $u \in \mc{S}_m$ and $K(u)= \inf_{\mc{S}_m} K$.
\end{itemize}
%
%
\end{Proposition}

%
%

\begin{Remark}
We highlight that with similar proofs we may show the same results of Propositions \ref{prop_uguag_not_sym}--\ref{prop_intr_L2min_PN}, \ref{prop_JJLU}, and Corollary \ref{corol_equival_min_min} in the radial setting, extending what we found in Theorem \ref{t-main} for $m>m_0$. 
\end{Remark}



We move now to some qualitative properties of the solutions. 
We highlight that, dealing with general solutions of a nonregular problem, their proofs are not standard and will require some additional arguments.
Since, in Theorem \ref{t-main}, we restricted to the radial framework, we investigate when this property is naturally achieved, together with positivity.


\begin{Proposition}[Symmetry and sign]\label{p-sym_unpert}
Assume \hyperref[g0]{\rm{(g0)}}--\hyperref[g2]{\rm{(g2)}} and let $\rho \in [-\infty,+\infty)$.
\begin{enumerate}[(i)]
\item All minimizers of $I^m$ over $\mc{P}_{(\rho)}$ 
and all minimizers of $K$ over $\mc{S}_m$ 
are radial about a point. 

\item All nonnegative minimizers of $I^m$ over $\mc{P}_{(\rho)}$ 
or $\mc{P}_{(\rho)}
^{\textup{rad}}$ 
and all nonnegative minimizers of $K$ over $\mc{S}_m$ 
or $\mc{S}_m
^{\textup{rad}}$ 
are Schwarz-symmetric up to a translation (see Definition \ref{def_schwartz}). 

\item Let $(\mu,u)$ be a solution to \eqref{eq_general_mu_fixed} 
such that $\intRN g_-(u)u \, \dx < +\infty$, and assume\footnote{Notice that if $\mu \geq 0$, this is implied by $g|_{(-\infty,0)} \ge 0$.} 
 $g(t)-\mu t \geq 0$ for each $t \leq 0$. 
 Then $u \geq 0$.


\item Assume, in addition, \hyperref[g3]{\rm{(g3)}} and \hyperref[g4]{\rm{(g4)}}. If $g|_{(-\infty,0)} = 0$ or $g$ is odd, then for any $m > m_0$ there exists a Schwarz-symmetric minimizer $(\mu_0, u_0)$ of $I^m$ over $\mc{P} 
$, where $\mu_0>0$ is also a minimizer of $\mu \in (-\infty,\overline{\mu}_0) \mapsto \inf_{\mc{P}_{\mu}} J_{\mu}-\frac{\mu}{2}m \in \R$ and $u_0$ is also a minimizer of $K$ over $\mc{S}_m$ and $\mc{D}_m$ with Lagrange multiplier $\mu_0$.
In particular,
\begin{align*}
 \inf_{ \mc{P} } I^m = d^m, \quad \inf_{ \mc{S}_m } K &= \kappa^m, \quad \inf_{\mc{D}_m} K = \ell^m, \\
 \inf_{-\infty< \mu < \overline{\mu}_0} \left( \inf_{\mc{P}_{\mu}} J_\mu - \mu \frac{m}{2}\right) =& \inf_{-\infty< \mu < \overline{\mu}_0} \left( a(\mu) - \mu \frac{m}{2}\right), 
\end{align*}
and all the above quantities are indeed equal and negative. 
This minimizer is a 
solution to \eqref{e-main} and satisfies the Nehari identity \eqref{eq_nehar_id}.
\end{enumerate}
\end{Proposition}

\begin{Remark}\label{r-ssm}
In a way analogous to the relation between Theorems \ref{t-main} and \ref{t-small-mass}, if either $g|_{(-\infty,0)} = 0$ or $g$ is odd, and  \hyperref[g0]{\rm{(g0)}}--\hyperref[g4]{\rm{(g4)}} and \eqref{eq_cond_meno_inf}\footnote{When $g|_{(-\infty,0)} = 0$, by \eqref{eq_cond_meno_inf} we mean $\lim_{t \to 0^+} g(t) / t = -\infty$. This convention will also be used in the proof of Proposition \ref{prop_esist_no_groundstate} below.} hold, then the outcome of Proposition \ref{p-sym_unpert} (iv), except for the signs of $\mu_0$ and $K(u_0)$ and the parts involving $\mc{D}_m$
, holds for every $m>0$.
\end{Remark}

When searching for a normalized solution without necessarily expecting minimality properties, we can drop \hyperref[g2]{\rm{(g2)}}; 
moreover, we can require the solution to be Schwarz symmetric (up to its sign).

\begin{Proposition}[Existence and symmetry under relaxed assumptions] 
\label{prop_esist_no_groundstate}
Assume \hyperref[g0]{\rm{(g0)}}, \hyperref[g1]{\rm{(g1)}}, 
\hyperref[g3]{\rm{(g3)}}, and \hyperref[g4]{\rm{(g4)}}.
Then, for every $m$ sufficiently large there exists a 
solution $(\mu_0, u_0) \in (0,+\infty) \times \mc{S}_m^\textup{rad}$ of \eqref{e-main} such that $K(u_0) < 0$ and $u_0$ or $-u_0$ is Schwarz symmetric
\footnote{$u_0$ if $t_0 > 0$, $-u_0$ if $t_0 < 0$, where $t_0$ is 
introduced in \hyperref[g4]{\rm{(g4)}}.}.
This solution satisfies the Poho\v{z}aev identity \eqref{eq_pohoz_frac} and the Nehari identity \eqref{eq_nehar_id}. Finally, if \eqref{eq_cond_meno_inf} holds as well, then for every $m>0$ there exists a solution $(\mu_0, u_0) \in \R \times \mc{S}_m^\textup{rad}$ of \eqref{e-main} that satisfies \eqref{eq_pohoz_frac} and \eqref{eq_nehar_id} and such that $u_0$ or $-u_0$ is Schwarz symmetric. 
\end{Proposition}

We highlight that in the final part of Proposition \ref{prop_esist_no_groundstate}, we do not have information on the sign of $\mu_0$. We end this section by giving some examples of nonlinearities to which the previous results 
apply. In these cases, we also show nonexistence of solutions.

\begin{Corollary}[Logarithmic case]\label{cor_log_power}
Let $m>0$ and
\[
g(t)= \alpha t \log(t^2) + \beta |t|^{q-1}t
\]
with $\alpha>0$, $\beta \in \R$, and $0 
< q \le 2^*_s - 1$. 
Then the following facts hold.
\begin{itemize}
	\item[(i)] If $q>1$ and $\beta \le -\frac{\alpha (q+1)}{q-1} e^{-(q+1)/2}$, then 
	there exist no solutions $(\mu,u) \in [0,+\infty) \times \mc{S}_m$ to \eqref{e-main} that satisfy the Poho\v{z}aev identity; in particular, no $L^2$-minima with nonnegative Lagrange multipliers $\mu$ exist. 
		\item[(ii)] If
		\begin{itemize}
			\item [(a)] $q\in (1, \bar{p})$ 
			and $\beta > 0$, or
			\item [(b)] $q \in (1, 2^*_s - 1]$ and $-\frac{\alpha (q+1)}{q-1} e^{-(q+1)/2} < \beta \le 0$, or
			\item [(c)] $q \in (0,1]$ and $\beta \leq 0$,
		\end{itemize}
		then there exists 
		$(\mu_0,u_0) \in \R 
		\times \mc{S}_m^\textup{rad}$, a 
		Schwarz-symmetric minimizer of $K$ over $\mc{S}_m$, which is a 
		solution of \eqref{e-main}, a Poho\v{z}aev minimum, and it 
		satisfies the Nehari identity. 
		Moreover, the following Legendre-transform formula holds:
		$$ 
		\kappa^m 
		= \inf_{-\infty< \mu < \overline{\mu}_0} \left( a(\mu) - \mu \frac{m}{2}\right) = a(\mu_0) - \mu_0 \frac{m}{2},
		$$ 
		where $\overline{\mu}_0 = \frac{\alpha}{q-1} \big[2\log\big(-\tfrac{\alpha}{\beta} \tfrac{q+1}{q-1}\big) - q - 1\big]$ in (b) with $\beta \ne 0$ and $\overline{\mu}_0 = +\infty$ otherwise. Finally, if $m_0$ is as in Theorem \ref{t-main} and $m>m_0$, then $\mu_0>0$ and $\ell^m=\kappa^m<0$.
\end{itemize}
\end{Corollary}

As already mentioned in Subsection \ref{subsec_literature}, we highlight that also the pure logarithmic case $\alpha=1$ and $\beta=0$, that is, \eqref{eq_intr_frac_log},
has an interest of its own: indeed, even if the existence of a normalized solution can be achieved by scaling arguments, in the fractional setting the existence of an $L^2$-minimum is not straightforward as in the local case, due to the lack of extremals for the fractional logarithmic Sobolev inequality.

\begin{Corollary}[Low-power case]
\label{corol_low_power}
Let
\[
g(t)= -\gamma |t|^{r-1} t + \beta |t|^{q-1}t
\]
with $\gamma, \beta>0$, $0 < r <1$, and $r < q < \bar{p}$.
Then, the conclusions of Corollary \ref{cor_log_power} (ii) hold. In particular, $\overline{\mu}_0 = +\infty$ if $q \in (1,\bar{p})$, $\overline{\mu}_0 = \beta
$ if $q = 1$, and $\overline{\mu}_0 = 2(q-r)\left(\frac{(1+r)(1-q)}{\gamma}\right)^{\frac{1-q}{q-r}} \left(\frac{\beta}{(1-r)(1+q)}\right)^{\frac{1-r}{q-r}}$ if $q \in (r,1)$.
\end{Corollary}

The results apply also to nonlinearities of the type $\alpha t \log(t^2) -\gamma |t|^{r-1} t + \beta |t|^{q-1}t$, but in this case the exact range of coefficients for existence is not explicitly available.

\begin{Remark}
	We notice that all the results we obtain are valid 
	in the local case $s=1$, $N\geq 3$, and this leads to some improvements of \cite{MeSc24}: in particular, in Theorem \ref{t-main}, the minimality on $\mc{P}^{\textup{rad}}$ of the solution, the Nehari identity, and formula \eqref{eq_legendre_intr} are new; moreover, Theorem \ref{t-small-mass} and all the results of Subsection \ref{subsec_furth_results} are new as well. 
	%
%
On the other hand, we do not cover the case $s=1$ and $N=2$, since some results would require some careful adaptations, mainly because $2^*_1 = \infty$ and the Poho\v{z}aev identity (cf. \cite[Remark 3.1]{BeLi83I} and \cite{BGK83}) leads to $\int G(u) - \frac{\mu}{2} u^2 \, \dx = 0$; see, in particular, the proofs of Lemma \ref{l-lambda-gen}, Lemma \ref{l-BJM} and Proposition \ref{p-a_monot} (notice that the existence of an $L^2$-minimum follows from \cite[Theorem 1.1]{MeSc24}).
However, other results contained here can be easily extended to this case; 
we leave the details to the interested reader.
%
\end{Remark}

\begin{Remark}
Some of the results contained in this paper 
are new even in the superlinear setting, extending or further developing \cite{CGT21N}. First, the growth of $g$ at infinity is 
slightly generalized (see \hyperref[g2]{\rm{(g2)}}-\hyperref[g3]{\rm{(g3)}}). Then, the techniques involved for the existence -- that is, the direct minimization over the Poho\v{z}aev set -- are new as well. Furthermore, new results appear, such as the fact that every Poho\v{z}aev minimum is a critical point (see Proposition \ref{pr-Poho_sol}), or the precise relation between $L^2$-minima and unconstrained ground states (see Proposition \ref{prop_JJLU}), together with other results in Subsection \ref{subsec_furth_results}. 
\end{Remark}







\paragraph{Outline of the paper.}
In Section \ref{s_prelim}, we prove some preliminary properties that will be used throughout the paper.
In Section \ref{s-BKS}, we recall and refine some results about the free-mass problem (that is, with fixed $\mu$).
In Section \ref{s-general}, together with Subsection \ref{subsec_symmetry_post}, we prove most of the results from Subsection \ref{subsec_furth_results}.
In Section \ref{s-pert}, we find a solution to the perturbed problem; then, in Section \ref{s:limit}, passing to the limit, we obtain a solution to \eqref{e-main} and prove Theorem \ref{t-main}.
Finally, in Section \ref{s-FER}, we prove Theorem \ref{t-small-mass}, Proposition \ref{prop_esist_no_groundstate}, and Corollaries \ref{cor_log_power} and \ref{corol_low_power}. 


\paragraph{Notations.} The definition of (weak) solution is given in Definition \ref{def_weak_sol}, while the definition of Schwarz-symmetric function is given in Definition \ref{def_schwartz}. The fractional Sobolev space $H^s(\R^N)$ is defined in \eqref{def_Ds_Hs}, while the restrictions to 
the radial subspace will be denoted with ``rad''. The $L^2$-sphere $\mc{S}_m$ and the $L^2$-ball $\mc{D}_m$ are given in \eqref{eq_def_Sm} and \eqref{eq_def_Dm_ellm}, while the $L^2$-critical exponent $\bar p + 1$ is defined in \eqref{eq_p}. 
The perturbed functions $g_{\eps}$ (and $\varphi_{\eps}$, $G^{\eps}_-$, $g^{\eps}_-$, $G_{\eps}$) are given in \eqref{eq_def_varphieps}--
\eqref{eq_def_geps}; the quantities related to the perturbed problem will appear with an ``$\eps$''.
The energy functional $K(u)$ and the $L^2$-minimal value (over the sphere or the ball) $\kappa^m$, $\ell^m$ are given in \eqref{eq_def_energ_func_sfer}, \eqref{eq_def_km_em}, \eqref{eq_def_Dm_ellm}.
The (action, Poho\v{z}aev) functionals $J_{\mu}(u)$, $P_{\mu}(u)$ over $H^s(\R^N)$ are defined in \eqref{eq_def_Jmu}, \eqref{eq_def_Pmu}, while
the (Lagrangian, Poho\v{z}aev) functionals $I^m(\mu, u)$, $P(\mu, u)$ over $\R \times H^s(\R^N)$
are given in \eqref{eq_def_Im}, \eqref{eq_def_Poh_pos}.
For technical reasons, the change of variable $\mu=e^{\lambda}$ in the frequency will be used when $\mu>0$, and sometimes an abuse of 
notation will be employed. In particular, we will write (leading to no ambiguity) 
$I^m(\lambda, u) \equiv I^m(e^{\lambda}, u)$, $P(\lambda, u) \equiv P(e^{\lambda}, u)$.
Thus, the (action, Lagrangian, Poho\v{z}aev) functionals $J(\lambda, u)$, $I^m(\lambda, u)$, $P(\lambda, u)$ over $\R \times H^s(\R^N)$ are given in \eqref{e-mpGamma} and at the beginning of Subsection \ref{subsec_minim_pohoz}.
The Poho\v{z}aev (set, minimal action) $\mc{P}_{\mu}$, $a(\mu)$ 
are given in \eqref{eq_def_poh_set}, \eqref{eq_def_amu}, 
while $\mc{P}$, $\mc{P}_{(\rho)}$, $\mc{P}_+$, $d^m$, $d^m_+$ 
can be found in \eqref{eq_def_mcP_pos}, \eqref{eq_def_pvarpi}, \eqref{eq_def_P+}, \eqref{eq_def_km_em}, and \eqref{eq_def_em+}. 
The frequency thresholds $\overline{\mu}_0$, $\overline{\lambda}_0$ appear in \eqref{eq_def_mubar0} and \eqref{e-mu0}.
The Legendre-type minimum $b^m$ is given in \eqref{e-infima}, while the 
mass threshold $m_0$ appears in \eqref{eq_meps_m0}. 
\section{Preliminaries}\label{s_prelim}


\indent

We first re-write here $g_{\pm}$ and $G_{\pm}$, introduced in \eqref{eq_intr_g+}-\eqref{eq_intr_g-}, in different terms: 
$$
G_+(t) :=
\begin{cases}
\displaystyle \int_0^t g^+(\tau) \, \mathrm{d}\tau & \text{if } t \ge 0,\\
\displaystyle \int_t^0 g^-(\tau) \, \mathrm{d}\tau & \text{if } t < 0,
\end{cases}
\quad
G_-(t):=
\begin{cases}
\displaystyle \int_0^t g^-(\tau) \, \mathrm{d}\tau & \text{if } t \ge 0,\\
\displaystyle \int_t^0 g^+(\tau) \, \mathrm{d}\tau & \text{if } t < 0,
\end{cases}
$$
$$g_+(t) := G_+'(t)
=
\begin{cases}
g^+ 
(t) & \text{if } t \ge 0,\\
-g^- 
(t) & \text{if } t \le 0,
\end{cases}
\quad g_-(t):= G_-'(t)
=
\begin{cases}
g^- 
(t) & \text{if } t \ge 0,\\
-g^+
(t) & \text{if } t \le 0.
\end{cases}
$$ 
Clearly, we have $G_{\pm}(t) \geq0$ and $g_{\pm}(t)t\geq 0$ for each $t \in\R$;
moreover,
$$G=G_+-G_-, \quad g=g_+-g_-,$$
and 
$$G_+(t) \geq G(t) \geq -G_-(t), \quad 
-g_-(t)t \le g(t)t \le g_+(t)t \quad \hbox{for $t \in \R$}.
$$
Generally, $G_{\pm}$ do not coincide with $G^{\pm}$, and 
unlike $G^{\pm}$, we have $G_{\pm} \in \mc{C}^1(\R)$. 
%
%
Finally, we notice that, if $g$ is odd, then $G_+, G_-$ are even.

\medskip

We briefly recall that the fractional Laplacian $(-\Delta)^s u := \mc{F}^{-1} (|\xi|^{2s} \mc{F}(u))$, $s \in (0,1)$, whenever $u$ is regular enough, can be expressed pointwise as
$$(-\Delta)^s u(x) = C_{N,s} \, \mathrm{P.V.} \int_{\R^N} \frac{u(x)-u(y)}{|x-y|^{N+2s}} \, \dy,$$
for some suitable $C_{N,s}>0$; we refer to \cite[Section 1.2]{Gal23}. 
We recall the fractional Sobolev space
\begin{equation}\label{def_Ds_Hs}
H^s(\R^N) := \left\{u \in L^2(\R^N) \mid |\cdot|^{s} \mc{F}(u) \in L^2(\R^N)\right\}
\end{equation}
and define
$$ \langle u,v \rangle :=\int_{\R^N} (-\Delta)^{s/2} u (-\Delta)^{s/2} v \, \dx = \frac{C_{N,s}}{2} \int_{\R^{2N}} \frac{\big(u(x)-u(y)\big)\big(v(x)-v(y)\big)}{|x-y|^{N+2s}} \, \dx \, \dy. $$
It is well-known that 
$H^s(\R^N) \hookrightarrow L^q(\R^N)$, $q \in [2, 2^*_s]$, while $H^s_{\textup{rad}}(\R^N) \hookrightarrow L^q(\R^N)$ is compact for $q \in (2, 2^*_s)$; here $H^s_{\textup{rad}}(\R^N)$ stands for the subspace of $H^s(\R^N) $ composed of radially symmetric functions.
We recall, moreover, the Gagliardo-Nirenberg inequality \cite{Par11GN} 
\begin{equation} \label{eq_GN_ineq}
|u|_{\bar{p}+1}^{\bar{p}+1} \leq C_\textup{GN} |(-\Delta)^{s/2} u|_2^2 |u|_2^{\bar{p}-1},
\end{equation}
where $C_\textup{GN}>0$ is optimal. 
We also recall some basic properties (see, e.g., 
\cite[Lemma 1.4.1]{Gal23}). 
\begin{Lemma}\label{l-assume+}
Let $u \in H^s(\R^N)$. 
Then $u^\pm,|u| \in H^s(\R^N)$ 
with $|(-\Delta)^{s/2}(u^{\pm})|_2 \leq |(-\Delta)^{s/2} u|_2$ and $|(-\Delta)^{s/2}|u||_2 \leq |(-\Delta)^{s/2} u|_2$. Finally, $\langle u,u^+ \rangle \ge \langle u^+ , u^+ \rangle$ and $\langle u,u^- \rangle \le - \langle u^- , u^- \rangle$.
\end{Lemma}


We show the following elementary convergence result.

\begin{Lemma}
\label{lem_converg_Hs}
Let $v \in H^s(\R^N)$. Let $v_R := v 
\varphi_R$, where $\varphi_R:= \varphi(\cdot/R)$ and $\varphi \in \mc{C}^{\infty}_c(\R^N)$ is chosen so that $0 \leq \varphi \leq 1$, $\varphi =1$ in $B_1(0)$, $\varphi=0$ in $\R^N \setminus B_2(0)$.
Then,
$$
\intRN (-\Delta)^{s/2} u (-\Delta)^{s/2} v_R \, \dx \to \intRN (-\Delta)^{s/2} u (-\Delta)^{s/2} v \, \dx \quad \text{as } R \to +\infty
$$
for each $u \in H^s(\R^N)$.
\end{Lemma}

\claim Proof.
We 
compute
\begin{align*}
\MoveEqLeft \frac{2
}{C_{N,s}} \intRN (-\Delta)^{s/2} u (-\Delta)^{s/2} (v_R-v) \, \dx\\
= & \, \int_{\R^{2N}} \frac{\big(u(x)-u(y)\big)\big((\varphi_R(x)-1) v(x) - (\varphi_R(y)-1) v(y)\big)}{|x-y|^{N+2s}} \, \dx \, \dy\\
= &\, \int_{\R^{2N}} (\varphi_R(x)-1)\frac{\big(u(x)-u(y)\big)\big( v(x) - v(y)\big)}{|x-y|^{N+2s}} \, \dx \, \dy\\
&+ \int_{\R^{2N}} \frac{\big(u(x)-u(y)\big)\big(\varphi_R(x) - \varphi_R(y)\big)}{|x-y|^{N+2s}} v(y) \, \dx \, \dy =: (I_1) + (I_2).
\end{align*}
As regards $(I_1)$, we observe that
\begin{equation*}
\frac{\big|u(x)-u(y)\big|\big| v(x) - v(y)\big|}{|x-y|^{N+2s}} =
\left( \frac{\big|u(x)-u(y)\big|^2}{|x-y|^{N+2s}} \right)^{1/2} \left( \frac{\big|v(x)-v(y)\big|^2}{|x-y|^{N+2s}} \right)^{1/2} \in L^1(\R^{2N}),
\end{equation*}
thus we can apply the dominated convergence theorem. Regarding $(I_2)$, we have
\begin{multline*}
\bigg | \int_{\R^{2N}} \frac{\big(u(x)-u(y)\big)\big(\varphi_R(x) - \varphi_R(y)\big)}{|x-y|^{N+2s}} v(y) \, \dx \, \dy \bigg | \\
\leq 
\left( \int_{\R^{2N}} \frac{\big|u(x)-u(y)\big|^2}{|x-y|^{N+2s}} \, \dx \, \dy \right)^{\frac{1}{2}} \left( \int_{\R^{2N}} \frac{\big|\varphi_R(x)-\varphi_R(y)\big|^2}{|x-y|^{N+2s} }v^2(y) \, \dx \, \dy \right)^{\frac{1}{2}},
\end{multline*}
and we can apply \cite[Lemma 1.4.5]{Amb21}.
%
\QED




\begin{Lemma}\label{lem_diff_along_L1}
Assume \hyperref[g0]{\rm{(g0)}}--\hyperref[g2]{\rm{(g2)}}. 
Then the functional
$$u \in H^s(\R^N) \mapsto \intRN G(u) \, \dx \in \R \cup \{-\infty\}$$
is differentiable at every $u$ 
such that $\intRN G(u) \, \dx > -\infty$ along every direction $\varphi \in H^s(\R^N) \cap L^1(\R^N)$.
\end{Lemma}

\claim Proof.
With the intent to apply Vitali's convergence theorem, we take $t \ne 0$ small and estimate 
\begin{align*}
\pabs{\frac{G(u+t\varphi)-G(u)}{t}} & = \pabs{\frac{1}{t} \int_{u}^{u+t \varphi} g(\tau) \, \mathrm{d}\tau}
\lesssim 
\pabs{\frac{1}{t} \int_{u}^{u+t \varphi} 1+ |\tau|^{2^*_s-1} \, \mathrm{d}\tau} \\
& \leq |\varphi| + \pabs{\frac{1}{t} \int_{u}^{u+t \varphi} |\tau|^{2^*_s-1} \, \mathrm{d}\tau}
\leq |\varphi| + \pabs{\frac{1}{t} \int_{u}^{u+t \varphi} (|u| + |\varphi|)^{2^*_s-1}\, \mathrm{d}\tau} \\
& = |\varphi| + (|u| + |\varphi|)^{2^*_s-1} |\varphi|
\lesssim |\varphi| + |u|^{2^*_s-1} |\varphi| + |\varphi|^{2^*_s},
\end{align*}
where we used the convexity of $t^{2^*_s-1}$. Then, fixed $A \subset \R^N$ measurable, we have
$$\int_A \pabs{ \frac{G(u+t\varphi)-G(u)}{t}} \, \dx \lesssim \int_{A} |\varphi| \, \dx + |u|_{2^*_s}^{2^*_s-1} \left(\int_{A} |\varphi|^{2^*_s} \, \dx\right)^{\frac{1}{2^*_s}} + \int_A |\varphi|^{2^*_s} \, \dx,$$
and thus 
we get both the uniform integrability and the tightness. This allows us to compute 
\begin{equation*}
\lim_{t \to 0} \frac{\intRN  G(u+t\varphi)\, \dx - \intRN G(u)\, \dx}{t} = \lim_{t \to 0} \intRN \frac{G(u+t\varphi)-G(u)}{t} \, \dx = \intRN g(u) \varphi \, \dx
\end{equation*}
and conclude.
%
\QED



\begin{Definition}
\label{def_weak_sol}
By \emph{weak solution} (or simply \emph{solution}) of \eqref{eq_general_mu_fixed} 
we mean $u\in H^s(\R^N)$ such that 
\begin{equation*}
\intRN (-\Delta)^{s/2} u (-\Delta)^{s/2} \varphi\, \dx + \mu \intRN u \varphi \, \dx = \intRN g(u) \varphi \, \dx
\end{equation*}
for each $\varphi \in H^s(\R^N) \cap L^1(\R^N)$; in particular, it holds for every $\varphi \in \mc{C}^{\infty}_c(\R^N)$.
\end{Definition}

We highlight that, generally, a 
solution $u$ of \eqref{eq_general_mu_fixed} is not a critical point of the action functional $J_{\mu}\colon H^s(\R^N) \to \R \cup \{+\infty\}$,
\begin{equation}\label{eq_def_Jmu}
J_{\mu}(u) := \frac12 \intRN |(-\Delta)^{s/2} u|^2 \, \dx + \frac{\mu}{2} \int_{\R^N} u^2 \, \dx - \intRN G(u) \, \dx,
\end{equation}
since, generally, $J_{\mu}$ is not differentiable at $u$ along every $\varphi \in H^s(\R^N)$ (see Lemma \ref{lem_diff_along_L1}); in particular, we cannot ensure that a 
solution satisfies the Nehari identity. On the other hand, in some circumstances, we will show that this is the case (cf. Proposition \ref{p-Nehari}).

Finally, we will use the following definition.
\begin{Definition}\label{def_schwartz}
A function $u \colon \R^N \to \R$ is called \emph{Schwarz-symmetric} if it is nonnegative, radially symmetric, and radially nonincreasing.
\end{Definition}

\section{The free-mass 
problem 
}\label{s-BKS}

\indent

The purpose of this section is to discuss the free-mass problem (i.e., $\mu$ fixed). 
For any $\mu \in \R$ we introduce a Poho\v{z}aev set
\begin{equation} \label{eq_def_poh_set}
\mc{P}_{\mu} := \left\{u \in H^s(\R^N) \setminus \{0\} \mid P_{\mu}(u) = 0\right\}, \quad \mc{P}_{\mu}^{\textup{rad}} := \mc{P}_{\mu} \cap H^s_{\textup{rad}}(\R^N),
\end{equation}
where the Poho\v{z}aev functional $
P_{\mu} \colon H^s(\R^N) \to \R \cup \{+\infty\}$ is defined by 
\begin{equation}\label{eq_def_Pmu}
P_{\mu}(u) 
:= \intRN |(-\Delta)^{s/2} u|^2 \, \dx + 2^*_s \intRN \frac{\mu}{2} u^2 - G(u) \, \dx.
\end{equation}
Recalled \eqref{eq_def_Jmu}, we 
define the least action
\begin{equation}\label{eq_def_amu}
	 a(\mu):= \inf_{\mc{P}_{\mu}^{\textup{rad}}} J_{\mu}.
	 \end{equation}
We notice that, generally, the set $\mc{P}_{\mu}$ could be empty, and in this case, $a(\mu)=+\infty$.

Here we list a series of properties that hold for $\mu \in \R$ 
under \hyperref[g0]{\rm{(g0)}}--\hyperref[g2]{\rm{(g2)}} (we recall from \eqref{eq_def_mubar0} the definition of $\overline{\mu}_0$). 

\begin{Lemma}\label{l-lambda-gen}
Assume \hyperref[g0]{\rm{(g0)}}--\hyperref[g2]{\rm{(g2)}} 
and let $\mu \in \R$. 
Then, the following statements are equivalent.
\begin{itemize}
	\item [(a)] $\mu < \overline{\mu}_0$.
	\item [(b)] There exists $t = t(\mu) \ne 0$ such that $G(t) >\frac{1}{2} \mu t^2 $.
	\item [(c)] There exists $u \in \mc{P}_{\mu}^{\textup{rad}}.$
\end{itemize}
\end{Lemma}

\claim Proof.
(a) $\iff$ (b) is trivial. 
If (c) holds, then $u \neq 
0$; thus, by the Poho\v{z}aev identity,
$$\int_{\R^N} 
G(u) - \frac{\mu}{2} u^2 
\, \dx = \frac{1}{2_s^*} \int_{\R^N} |(-\Delta)^{s/2} u|^2 \, \dx > 0,$$
hence there exists $t(\mu):=u(x_{\mu})$, $x_{\mu} \in \R^N$, such that (b) is verified.
Assume (b) holds, from \cite[Proof of Theorem 2]{BeLi83I} (see also \cite[Lemma 3.1]{CGT21N}) 
there exists $\bar{u} \in H^1_{\textup{rad}}(\R^N) \subset H^s_{\textup{rad}}(\R^N)$ such that $\intRN \big(G(\bar{u}) - \frac{\mu}{2} \bar{u}^2 \big)\, \dx > 0$. Then, for $t>0$, straightforward computations show that
\[
P_{\mu}\bigl(\bar{u}(\cdot/t)\bigr) =\left( \intRN |(-\Delta)^{s/2} \bar u|^2 \, \dx \right) t^{N-2s} + 2^*_s \left(\intRN \frac{\mu}{2} \bar u^2 - G(\bar u) \, \dx \right) t^N
\]
is positive when $t$ is small and negative when $t$ is large, hence there exists $\bar{t}> 0$ such that $P_{\mu}\bigl(\bar{u}(\cdot/\bar{t})\bigr) = 0$, whence (c). 
%
\QED

\begin{Remark}\label{r-all_mu}
By the previous proof we notice the following fact: 
for every $u \in H^s(\R^N) \setminus \{0\}$ such that $P_\mu(u) \le 0$, we have $\intRN \big( G(u) - \frac{\mu}{2} u^2\big) \, \dx > 0$, thus there exists $t_0>0$ such that $u_0 := u(\cdot/t_0) \in \mc{P}_\mu$.
Moreover, using that $P_{\mu}(u) \le 0$ and that $t \in (0,\infty) \mapsto \frac{1}{2} t^{N-2s}- \frac{1}{2^*_s} t^N\in \R$ is maximized at $t=1$, there holds
\begin{align*}
\frac{s t_0^{N-2s}}{N} \intRN |(-\Delta)^{s/2} u|^2 \, \dx & = \frac{s}{N} \intRN |(-\Delta)^{s/2} u_0|^2 \, \dx = J_{\mu}(u_0) \\
& = \frac{t_0^{N-2s}}{2} \intRN |(-\Delta)^{s/2} u|^2 \, \dx - t_0^N \intRN G(u) - \frac{\mu}{2} u^2 \, \dx \\
& \le \left(\frac{t_0^{N-2s}}{2} - \frac{t_0^N}{2^*_s}\right) \intRN |(-\Delta)^{s/2} u|^2 \, \dx \le \frac{s}{N} \intRN |(-\Delta)^{s/2} u|^2 \, \dx,
\end{align*}
that is, $\intRN |(-\Delta)^{s/2} u_0|^2 \, \dx \le \intRN |(-\Delta)^{s/2} u|^2 \, \dx$, 
with a strict inequality if $P_{\mu}(u) < 0$, and $t_0 \le 1$.
\\
\noindent Notice that this fact does not depend on the sign of $\mu$.
\end{Remark}

\begin{Lemma}\label{l-BJM}
Let $\mu \in \R$ and assume \hyperref[g0]{\rm{(g0)}}--\hyperref[g2]{\rm{(g2)}}. 
If $\mc{P}_\mu\ne \emptyset$ and $u \in \mc{P}_\mu$ is such that $J_\mu(u) = \inf_{\mc{P}_{\mu}} J_\mu$, then $\omega := \intRN G(u) - \frac{\mu}{2} u^2 \, \dx > 0$ and
\begin{equation}\label{eq_minBS}
\frac12 \abs{(-\Delta)^{s/2}u}_2^2 = \min\left\{\frac12 \abs{(-\Delta)^{s/2}v}_2^2 \mid v \in H^s(\R^N), \; \intRN G(v) - \frac{\mu}{2} v^2 \, \dx = \omega \right\}.
\end{equation}
\end{Lemma}
\claim Proof.
First, notice that since $\mc{P}_\mu\ne \emptyset$, there exists a function $v\in H^s(\R^N)$ such that $\intRN G(v) - \frac{\mu}{2} v^2 \, \dx = 1$, thus the minimization problem 
\[
T := \inf\left\{\frac12 \abs{(-\Delta)^{s/2}v}_2^2 \mid v \in H^s(\R^N), \; \intRN G(v) - \frac{\mu}{2} v^2 \, \dx = 1 \right\}
\]
is well defined. Moreover, every such function verifies $G(v)\in L^1(\R^N)$. 
We adapt now the proof of \cite[Lemma 1 (ii)]{BJM09}. We notice that we do not need to assume that the Poho\v{z}aev identity holds for every solution or that $T$ is a priori attained.
Since $u \in \mc{P}_\mu$, we have that
\[
J_\mu(u) \ge 2s (N-2s)^{\frac{N}{2s}-1} N^{-\frac{N}{2s}} T^\frac{N}{2s} =: \widetilde{T}.
\]
Now, let $(v_n)_n$ be a minimizing sequence for $T$, and let $t_n > 0$ be such that $v_n(\cdot/t_n) \in \mc{P}_\mu$. Then, since $J_\mu(u) = \inf_{\mc{P}_{\mu}} J_\mu$,
$$
\widetilde{T} \le J_\mu(u) \le J_\mu\bigl(v_n(\cdot/t_n)\bigr) = \widetilde{T} + o_n(1),
$$
and we conclude letting $n \to +\infty$ and via scaling arguments.
\QED


\begin{Proposition}\label{p-a_monot}
Assume \hyperref[g0]{\rm{(g0)}}--\hyperref[g2]{\rm{(g2)}} and let $\mu \in \R$. 
Then $a(\mu) \geq 0$,
\begin{equation}\label{eq_min_poh_less}
a(\mu) = \inf \bigg\{\frac{s}{N} \intRN |(-\Delta)^{s/2} u|^2 \, \dx \mid u \in H^s_{\textup{rad}}(\R^N) \setminus\{0\}, \; P_{\mu}(u)\leq 0 \bigg\},
\end{equation}
and $a(\mu) < \frac{s}{N} \intRN |(-\Delta)^{s/2} u|^2 \, \dx$ if $P_{\mu}(u) < 0$. 
Moreover, $a \colon \R \to [0,+\infty]$ is nondecreasing. Finally, if $\mu > 0$ and $\mc{P}^\textup{rad}_\mu \ne \emptyset$, then $a$ is continuous at $\mu$.
\end{Proposition}

When $s=1$, \eqref{eq_min_poh_less} was already noted in \cite[Proposition 1.5]{Sha85}.

\smallskip

\claim Proof of Proposition \ref{p-a_monot}.
By definition of $\mc{P}^\textup{rad}_{\mu}$, we have
$$a(\mu) = \frac{s}{N} \inf_{u \in \mc{P}_{\mu}^{\textup{rad}}} \abs{(-\Delta)^{s/2} u}_2^2 \geq 0.$$
To show relation \eqref{eq_min_poh_less}, observe that if $\mc{P}^\textup{rad}_\mu = \emptyset$, then $\bigl\{u \in H^s_\textup{rad}(\R^N) \setminus \{0\} \mid P_\mu(u) \le 0 \bigr\} = \emptyset$ as well from Remark \ref{r-all_mu}, hence both sides of \eqref{eq_min_poh_less} equal $+\infty$. Next, assume $P_{\mu}(u) \leq 0$; thus, by Remark \ref{r-all_mu} there exists $t_0>0$ such that $P_{\mu}\bigl(u(\cdot/t_0)\bigr) = 0$ and, consequently,
$$a(\mu) \leq \frac{s}{N} \intRN |(-\Delta)^{s/2} u(\cdot/t_0)|^2 \, \dx \leq \frac{s}{N} \intRN |(-\Delta)^{s/2} u|^2 \, \dx$$
(where the second inequality is strict if $P_{\mu}(u) < 0$), 
and we can pass to the infimum. The reverse inequality is obvious. 
We show now the monotonicity. Let $\mu_1 < \mu_2$. If $\mc{P}^\textup{rad}_{\mu_2}=\emptyset$, then it clearly holds. Otherwise, let $u \in \mc{P}^\textup{rad}_{\mu_2}$. Then $P_{\mu_1}(
u) < P_{\mu_2}(u) = 0$ and by \eqref{eq_min_poh_less} we have 
\begin{equation}\label{eq_dim_monot_amu}
a(\mu_1) \leq \frac{s}{N} \intRN |(-\Delta)^{s/2} u|^2 \, \dx = J_{\mu_2}(u);
\end{equation}
passing to the infimum, we have the claim. 
We show now the continuity, starting with the upper semicontinuity. 
Let $\mu_n \to \mu>0$, $\eps>0$ fixed, and $u \in \mc{P}^\textup{rad}_{\mu}$ such that $J_{\mu}(u) \leq a(\mu)+\eps$. From Remark \ref{r-all_mu} we have $\int_{\R^N} \left(G(u)-\frac{\mu}{2} u^2\right)\, \dx>0$ and thus, for $n$ large, $\int_{\R^N} \left(G(u)-\frac{\mu_n}{2} u^2\right)\, \dx>0$, which means that there exists $t_n>0$ such that $u_n:=u(\cdot/t_n) \in \mc{P}_{\mu_n}$; notice that $t_n \to 1$. Hence,
$$a(\mu_n) \leq J_{\mu_n}(u_n) = \frac{s}{N} \norm{(-\Delta)^{s/2} u_n}_2^2 = \frac{s}{N} t_n^{N-2s} \norm{(-\Delta)^{s/2} u}_2^2 = t_n^{N-2s} J_{\mu}(u) \leq t_n^{N-2s} \big( a(\mu) +\eps\big),$$
and the claim follows by letting $n \to +\infty$ and the arbitrariness of $\eps$. Now, 
we focus on the lower semicontinuity and fix $\mu_n \to \mu>0$. From \hyperref[g0]{\rm{(g0)}} and Lemma \ref{l-lambda-gen}, $\mc{P}^\textup{rad}_{\mu_n} \ne \emptyset$ for all sufficiently large $n$; in addition, we can assume $\mu^*:= \inf_n \mu_n >0$, 
and by upper semicontinuity, $M:=\sup_n a(\mu_n) <\infty$.
Fix $\eps > 0$. For every $n$ there exists $u_n \in \mc{P}^\textup{rad}_{\mu_n}$ such that $J_{\mu_n}(u_n) \le a(\mu_n) + \eps$.
We show that $(u_n)_n$ is bounded. 
Indeed, by the Poho\v{z}aev identity we have
$$\frac{s}{N}|(-\Delta)^{s/2} u_n|_2^2 = J_{\mu_n}(u_n) \leq M + \eps$$
and, again by the Poho\v{z}aev identity and the assumptions about $G_+$,
$$ \frac{\mu^*}{2} |u_n|_2^2 \le \frac{\mu_n}{2} |u_n|_2^2 \leq \int_{\R^N} G(u_n) \, \dx + C \leq  \frac{\delta}{2} |u_n|_2^2+ C_{\delta} |u_n|_{2^*_s}^{2^*_s} 
+C \leq\frac{\delta}{2} |u_n|_2^2+ C,$$
where, by choosing $\delta<\mu^*$, we obtain the boundedness of $u_n$ in $H^s(\R^N)$. 
From Remark \ref{r-all_mu} we have $
\int_{\R^N} \left(G(u
_n) - \frac{\mu_n}{2} u
_n^2\right)\, \dx = \frac{1}{2^*_s}|(-\Delta)^{s/2}u_n|_2^2 >0$; moreover, it is bounded away from zero, since (again choosing $\delta<\mu^*$)
$$ \frac{1}{2^*_s}|(-\Delta)^{s/2}u_n|_2^2 = - \frac{\mu_n}{2} |u_n|_2^2 +\int_{\R^N} G(u_n) \, \dx \leq \frac{\delta-\mu_n}{2} |u_n|_2^2 + C_{\delta} |u_n|_{2^*_s}^{2^*_s} \leq C_{\delta}' |(-\Delta)^{s/2}u_n|_2^{2^*_s}.$$
This and the boundedness of $u_n$ in $L^2(\R^N)$ yield $\int_{\R^N} \left(G(u_n) - \frac{\mu}{2} u_n^2\right)\, \dx>0$ for $n$ large, hence there exists $t_n >0$ such that $v_n:=u_n(\cdot/t_n) \in \mc{P}_{\mu}$, which implies
$$
a(\mu) \leq J_{\mu}(v_n) = \frac{s}{N} |(-\Delta)^{s/2} v_n|_2^2 = \frac{s}{N} t_n^{N-2s} |(-\Delta)^{s/2} u_n|_2^2 = t_n^{N-2s} J_{\mu_n}(u_n) \leq t_n^{N-2s} \left(a(\mu_n) + \eps\right)
$$
and the claim follows by letting $n \to +\infty$ (noticed that $t_n \to 1$) and the arbitrariness of $\eps$. 
\QED

%



\begin{Proposition}\label{pr_bks}
Assume \hyperref[g0]{\rm{(g0)}}--\hyperref[g2]{\rm{(g2)}} and let $\mu \in \R$.
\begin{enumerate}[(i)]
	\item All minimizers of $J_{\mu}$ over $\mc{P}_{\mu}$ or $\mc{P}_{\mu}^{\textup{rad}}$ are 
	solutions to
	\begin{equation}\label{e-massfree}
		(-\Delta)^s u + \mu u = g(u) \quad \text{in } \R^N.
	\end{equation}
	
	\item All minimizers of $J_{\mu}$ over $\mc{P}_{\mu}$ are radially symmetric about a point.
	
	\item All nonnegative minimizers of $J_{\mu}$ over $\mc{P}_{\mu}$ or $\mc{P}_{\mu}^{\textup{rad}}$ are Schwarz-symmetric up to a translation.
	
	\item If $g(t) - \mu t \ge 0$ for all $t \le 0$, then every 
	solution to \eqref{e-massfree} such that $\intRN g_-(u) u \, \dx < +\infty$ is nonnegative. 
\end{enumerate}
\end{Proposition}

\begin{Remark}
In Proposition \ref{pr_bks} (i), we are not claiming that, for general $g$, there exists a minimizer over $\mc{P}_{\mu}$; on the other hand, from Proposition \ref{pr_bks} (ii), we deduce that $\mc{P}_{\mu}^{\textup{rad}}$ is a natural subset where to search for such a minimizer.
\end{Remark}

The proof of Proposition \ref{pr_bks} is analogous to those of Propositions \ref{pr-Poho_sol} and \ref{p-sym_unpert} (see Subsections \ref{subsub_l2_poh_crit}, \ref{subsec_symmetry_post} and \ref{subsec_passage_lim}), hence it is omitted.
We now state the main theorem of this section.

\begin{Theorem}\label{t-bks}
Let $\mu > 0$ and $g \in \mc{C}(\R)$ satisfy
\begin{itemize}
\item
$|g(t)| \lesssim |t| + |t|^{2^*_s-1}$ for every $t \in \R$;
\item
$\displaystyle \limsup_{t \to 0} \frac{G(t)}{t^2}\leq 0$ and $\displaystyle \limsup_{|t| \to \infty} \frac{G(t)}{|t|^{2^*_s}} \leq 0$;
\item
there exists $
t_0 \ne 0$ such that $G(t_0) > \frac{1}{2}\mu t_0^2$.
\end{itemize}
Then the following statements hold true.
\begin{enumerate}[(i)]
\item There exists $v_0 \in \mc{P}_{\mu}^{\textup{rad}}$ such that 
\begin{equation*}
J_{\mu}(v_0) =
a(\mu)
\end{equation*}
and thus $a(\mu)> 0$.
%
%
Moreover, $a \colon (0,\overline{\mu}_0) \to (0,+\infty)$ is increasing.

\item If $g|_{(-\infty,0)} = 0$ or $g$ is odd, then there exists a Schwarz-symmetric minimizer of $J_{\mu}$ over $\mc{P}_{\mu}$; in particular,
\[
\inf_{\mc{P}_{\mu}} J_{\mu} = \inf_{\mc{P}_{\mu}^{\textup{rad}}} J_{\mu}.
\]
%
\end{enumerate}
\end{Theorem}

\claim Proof.

(i) 
%
From the assumptions, we see that for every $\delta > 0$ there exists $C_\delta > 0$ such that for every $t \in \R$
\begin{equation*}
	G(t) \le \delta t^2 + C_\delta |t|^{2^*_s},
\end{equation*}
whence, for every $u \in H^s(\R^N)$,
\begin{equation*}
\frac{1}{2^*_s} P_{\mu}(u) \ge \frac{1}{2^*_s} |(-\Delta)^{s/2} u|_2^2 + \frac\mu2 |u|_2^2 - \delta |u|_2^2 - C_\delta |u|_{2^*_s}^{2^*_s} \ge \left(\min\left\{\tfrac{1}{2^*_s},\tfrac{\mu}{2}\right\} - \delta\right) \norm{u}_{H^s}^2 - C_\delta \norm{u}_{H^s}^{2^*_s}.
\end{equation*}
Taking $0 < \delta < \min\{1/2^*_s,\mu/2\}$, we see that $P_{\mu}(u) > 0$ if $\norm{u}_{H^s}$ is sufficiently small.

Let now $(u_n)_n \subset \mc{P}_{\mu}^{\textup{rad}}$ such that $\lim_n J_{\mu}(u_n) = \inf_{\mc{P}_{\mu}^{\textup{rad}}} J_{\mu}$. Since $P_{\mu}(u_n) = 0$, we see that
\begin{equation*}
	\inf_{\mc{P}_{\mu}^{\textup{rad}}} J_{\mu} + o(1) = J_{\mu}(u_n) - \frac{1}{2^*_s} P_{\mu}(u_n) = \frac{s}{N} |(-\Delta)^{s/2} u_n|_2^2,
\end{equation*}
which yields that $|(-\Delta)^{s/2} u_n|_2$ is bounded, and, using the same $\delta$-related argument as before, that
\begin{align*}
	\inf_{\mc{P}_{\mu}^{\textup{rad}}} J_{\mu} + o(1) & = J_{\mu}(u_n) - \frac{1}{2} P_{\mu}(u_n) = -\frac{\mu}{2} \left(\frac{2^*_s}{2} - 1\right) |u_n|_2^2 + \left(\frac{2^*_s}{2} - 1\right) \intRN G(u_n) \, \dx\\
	& \le -\frac{\mu}{2} \left(\frac{2^*_s}{2} - 1\right) |u_n|_2^2 + \left(\frac{2^*_s}{2} - 1\right) \left(\delta |u_n|_2^2 + C_\delta |u_n|_{2^*_s}^{2^*_s}\right).
\end{align*}
Since $|u_n|_{2^*_s}$ is bounded from the previous step and the fractional Sobolev embedding, taking $\delta < \mu/2$ we obtain that $|u_n|_2$ is bounded as well, i.e., $\norm{u_n}_{H^s}$ is bounded. 
We can then consider a subsequence, still denoted by $u_n$, and $u_0 \in H_\textup{rad}^s(\R^N)$ such that $u_n \rightharpoonup u_0$ in $H^s(\R^N)$, $u_n \to u_0$ a.e. in $\R^N$, and $u_n \to u_0$ in $L^q(\R^N)$ for every $q \in (2,2^*_s)$. In particular, using the weak lower semi-continuity of the norm and Fatou's lemma, we get that $J_{\mu}(u_0) \le \inf_{\mc{P}_{\mu}^{\textup{rad}}} J_{\mu}$ and $P_{\mu}(u_0) \le 0$. In addition, recalling $P_{\mu}(u_n) = 0$, we have by the strong convergence in $L^{q}(\R^N)$, 
\begin{align*}
\MoveEqLeft	|(-\Delta)^{s/2} u_0|_2^2 + \frac{2^*_s \mu}{2} |u_0|_2^2 + 2^*_s \intRN G_-(u_0) \, \dx\\
	&\le \limsup_n\left( |(-\Delta)^{s/2} u_n|_2^2 + \frac{2^*_s \mu}{2} |u_n|_2^2 + 2^*_s \intRN G_-(u_n) \, \dx \right)\\
	&= 2^*_s \limsup_n \intRN G_+(u_n) \, \dx = 2^*_s \intRN G_+(u_0) \, \dx; 
\end{align*}
if $u_0 \equiv 0$, the above relation implies 
$u_n \to 0$ in $H^s(\R^N)$, which is impossible because $P_{\mu}(u) > 0$ for small $u$. Since $u_0 \ne 0$, by Remark \ref{r-all_mu} 
we can find $t_0>0$ such that $P_{\mu}\bigl(u_0(\cdot/t_0)\bigr) = 0$ and
\begin{align*}
J_{\mu}\bigl(u_0(\cdot/t_0)\bigr) & = \frac{s}{N} \intRN |(-\Delta)^{s/2} u_0(\cdot/t_0)|^2 \, \dx \le \frac{s}{N} \intRN |(-\Delta)^{s/2} u_0|^2 \, \dx\\
& \le \lim_n \frac{s}{N} \intRN |(-\Delta)^{s/2} u_n|^2 \, \dx = \lim_n J_\mu(u_n) = a(\mu),
\end{align*}
so the sought minimizer is $v_0 := u_0(\cdot/t_0)$.
By the Poho\v{z}aev identity we also have $a(\mu) = J_{\mu}(v_0) = \frac{s}{N} |(-\Delta)^{s/2} v_0|_2^2 > 0$.

To show the strict monotonicity, arguing as in \eqref{eq_dim_monot_amu}, in the case $\mu_2>0$ we choose $u$ as the minimum corresponding to $\mu_2$. Hence, 
from Proposition \ref{p-a_monot}, $P_{\mu_1}(u) < P_{\mu_2}(u) = 0$ and 
\[
a(\mu_1) < \frac{s}{N} \intRN |(-\Delta)^{s/2} u|^2 \, \dx = J_{\mu_2}(u) = a(\mu_2).
\]

\smallskip

(ii) We start with the case when $g$ is odd. Let $(u_n)_n \subset \mc{P}_{\mu}$ be such that $\lim_n J_{\mu}(u_n) = \inf_{\mc{P}_{\mu}} J_{\mu}$. From Lemma \ref{l-assume+}, $J_{\mu}(|u_n|) \le J_{\mu}(u_n)$ and $P_{\mu}(|u_n|) \le P_{\mu}(u_n) = 0$. If $t_n \in (0,1]$ 
is such that $P_{\mu}\bigl(|u_n|(\cdot/t_n)\bigr) = 0$, cf. Remark \ref{r-all_mu}, then we obtain 
$$\inf_{\mc{P}_{\mu}}J_{\mu} \le J_{\mu}\bigl(|u_n|(\cdot/t_n)\bigr) = \frac{s}{N} t_n^{N-2s} |(-\Delta)^{s/2}|u_n||_2^2 \leq   \frac{s}{N}|(-\Delta)^{s/2} u_n|_2^2= J_{\mu}(u_n).$$
This proves that, up to replacing $u_n$ with $|u_n|(\cdot/t_n)$, we can assume that $u_n \ge 0$ for every $n$. In a very similar way, using the Schwarz rearrangement and the fractional 
P\'olya--Szeg\"o inequality \cite{Par11PS}, 
we can assume that 
every $u_n$ is Schwarz symmetric. We can now conclude in the very same way as in point (i). 
If $g|_{(-\infty,0)} = 0$, we argue as above, except we use $u_n^+$ instead of $|u_n|$.
\QED

\begin{Remark}
Theorem \ref{t-bks} (i) is a refinement of \cite[Theorem 1.2 (i)]{BKS17}, where the main additional assumption is the existence of $q \in (2,2^*_s)$ such that $\lim_{t \to 0} g(t) / |t|^{q-1} = 0$. At the same time, Proposition \ref{pr_bks} shows that \cite[Theorem 1.2 (ii)--(iii)]{BKS17} holds under considerably weaker assumptions (in particular, \hyperref[g1]{\rm{(g1)}}).
\end{Remark}

\begin{Remark}
Assume \hyperref[g0]{\rm{(g0)}}--\hyperref[g2]{\rm{(g2)}} and that $\mu_0$ is a \emph{global minimum} of the function
$\mu \mapsto  a(\mu) - \frac{\mu}{2} m$.
Then, the existence of an action minimum for $a(\mu_0)$
(as given, for instance, in Theorem \ref{t-bks}), implies the existence of a (radial) $L^2$-minimum with mass $m$. 
Indeed, 
exploiting 
Propositions \ref{prop_uguag_not_sym} and \ref{prop_intr_legendre}, 
the minimality of $\mu_0$
and of $u_0$,
the definitions of $I^m$ and $d^m$, and 
the fact that $(\mu_0, u_0) \in \mc{P}$, 
we obtain
\begin{equation*}
\kappa^m = \inf_{\mu < \overline{\mu}_0} \left( a(\mu) - \frac{\mu}{2} m \right) = a(\mu_0) - \frac{\mu_0}{2} m = J_{\mu_0}(u_0) - \frac{\mu_0}{2} m = I^m(\mu_0,u_0) \ge d^m = \kappa^m,
\end{equation*}
thus $I^m(\mu_0,u_0) = d^m$ and hence $(\mu_0, u_0)$ is a Poho\v{z}aev minimum, which implies, by Corollary \ref{corol_equival_min_min},  that $u_0$ is an $L^2$-minimum. 

More generally, one may wonder if a \emph{critical point} of $\mu \mapsto  a(\mu) - \frac{\mu}{2} m$ may lead to a normalized solution with mass $m$. This nontrivial question has been tackled recently in the local framework (see \cite{CGT26} and references therein), by considering local minima and local maxima of such a map. 
A similar argument, albeit in a mildly different context, is used in \cite{DC_D_G_S}. 

\end{Remark}

\section{Some general properties}\label{s-general}

\indent

For $\rho \in [-\infty,+\infty)$, recall $\mc{P}_{(\rho)}$ from \eqref{eq_def_pvarpi} and $\overline{\mu}_0$ from \eqref{eq_def_mubar0}. 
%
%
We observe that, by Lemma \ref{l-lambda-gen}, if $\rho \in [-\infty, \overline{\mu}_0)$, then $\mc{P}_{(\rho)} \supset \mc{P}_{(\rho)}^{\textup{rad}}
\ne \emptyset$; in particular, $\mc{P} 
\ne \emptyset$ always holds.

\begin{Remark}\label{rem_2min}
We point out that if $(\mu,u)$ minimizes $I^m$ over $\mc{P}_{(\rho)}$ for some $\rho \in [-\infty,+\infty)$, then $u$ minimizes $J_\mu$ over $\mc{P}_\mu$.
\end{Remark}

\subsection{Connection among infima}

\indent

We prove now Propositions \ref{prop_uguag_not_sym} and \ref{prop_intr_legendre}.

\medskip

\claim Proof of Proposition \ref{prop_uguag_not_sym}.
Let $(\mu_n, u_n)_n \subset \mc{P}$ be a minimizing sequence for $\inf_{\mc{P}} I^m$, that is, $I^m(\mu_n, u_n) =\inf_{\mc{P}} I^m + o(1)$; 
recall that $u_n \neq 0$ by the definition of $\mc{P}$. Consider $t_n:= (m/|u_n|_2^2)^{1/N}$ such that $v_n:=u_n(\cdot / t_n) \in \mc{S}_m$. We notice, by the Poho\v{z}aev identity,
\begin{align*}
I^m(\mu_n, u_n) &= \frac{1}{2} |(-\Delta)^{s/2}u_n|_2^2 + \frac{\mu_n |u_n|_2^2}{2} \left(1 - t_n^N\right) -\int_{\R^N}G(u_n) \, \dx \\
& =\frac{1}{2} |(-\Delta)^{s/2}u_n|_2^2 + \left(1 - t_n^N\right) \left( \int_{\R^N} G(u_n) \, \dx - \frac{1}{2^*_s} |(-\Delta)^{s/2} u_n|_2^2\right)-\int_{\R^N}G(u_n) \, \dx \\
& =\left(\frac{1}{2} - \frac{1}{2^*_s} \left(1 - t_n^N\right)\right) |(-\Delta)^{s/2}u_n|_2^2 - t_n^N \int_{\R^N} G(u_n) \, \dx \\
& =\left(\frac{1}{2} - \frac{1}{2^*_s} \left(1 - t_n^N\right)\right) \frac{1}{t_n^{N-2s}} |(-\Delta)^{s/2} v_n|_2^2 - \int_{\R^N} G(v_n) \, \dx \\
& \geq \frac12 |(-\Delta)^{s/2} v_n|_2^2 - \int_{\R^N} G(v_n) \, \dx = K(v_n),
\end{align*}
where we used that $\frac{1}{2} - \frac{1}{2^*_s} \left(1 - t^N\right) \geq \frac12 t^{N-2s}$ for all $t>0$, which is true by a straightforward computation.
Thus, we have shown
$$
\inf_{\mc{S}_m} K \leq K(v_n) \leq I^m(\mu_n, u_n) = \inf_{\mc{P}} I^m + o(1).$$
Passing to the limit in $n$, we obtain
\begin{equation}\label{eq_KI_in}
\inf_{\mc{S}_m} K \leq \inf_{\mc{P}} I^m.
\end{equation}

Assume now $\rho \in [0, +\infty)$ and $\inf_{\mc{D}_m}K<-\frac{\rho}{2} m$. We can thus choose a minimizing sequence $(u_n)_n \subset \mc{D}_m$ such that $\lim_n K(u_n) = \inf_{\mc{D}_{m}} K$ and $K(u_n)<-\frac{\rho}{2} m \le 0$ for every $n$ (which implies $u_n \ne 0$)
and define
\begin{equation*}
\mu_n := \frac{2}{N |u_n|_2^2} \left(s|(-\Delta)^{s/2} u_n|_2^2 - N K(u_n)\right) \geq 
-\frac{2}{ |u_n|_2^2} K(u_n) > \rho
\end{equation*}
so that $(\mu_n,u_n) \in \mc{P}_{(\rho)}$.
This, together with $\mu_n >0$ and $|u_n|_2^2 \le m$, yields
\[
\inf_{\mc{D}_{m}} K + o_n(1) = K(u_n) \ge I^m(\mu_n,u_n) \ge \inf_{\mc{P}_{(\rho)}} I^m,
\]
and passing to the limit in $n$ we get $\inf_{\mc{D}_{m}} K \ge \inf_{\mc{P}_{(\rho)}} I^m$.
Arguing similarly for $\rho \in 
\R%
$ and $\inf_{\mc{S}_m}K<-\frac{\rho}{2} m$ we obtain 
$\inf_{\mc{S}_{m}} K \ge \inf_{\mc{P}_{(\rho)}} I^m.$
Choosing $\rho=-\infty$ and recalling \eqref{eq_KI_in}, we get \eqref{eq_ineq_not_sym}.
Gathering the previous inequalities, we conclude the proof.

Finally, we show the last relation arguing as in the proof of \eqref{eq_min_poh_less}. Indeed, considered $(\mu, u) \in (\rho, +\infty) \times H^s(\R^N)$ with $P(\mu, u) \leq 0$, there exists $t_0 > 0$ such that $P(\mu, u(\cdot/t_0))=0$ and $|(-\Delta)^{s/2} u(\cdot/t_0)|_2^2 \leq |(-\Delta)^{s/2} u|_2^2$ 
(with a strict inequality if $P(\mu, u) < 0$), 
thus
$$\inf_{\mc{P}_{(\rho)}} I^m \leq I^m(\mu, u(\cdot/t_0)) = \frac{s}{N} |(-\Delta)^{s/2} u(\cdot/t_0)|_2^2 -\frac{\mu}{2} m \leq \frac{s}{N} |(-\Delta)^{s/2} u|_2^2 -\frac{\mu}{2} m$$
and hence the claim passing to the infimum (the other inequality is obvious).
\QED


\bigskip


\medskip

\claim Proof of Proposition \ref{prop_intr_legendre}.
Let $\mu \in (\rho, \overline{\mu}_0)$. From Lemma \ref{l-lambda-gen}, $\mc{P}_{\mu} \neq \emptyset$ (thus $\inf_{\mc{P}_{\mu}} J_{\mu}$ is well defined).
Then, for any $u \in \mc{P}_{\mu}$, we have
$$\inf_{\mc{P}_{(\rho)}} I^m \leq I^m(\mu, u) = J_{\mu}(u) - \frac{\mu}{2} m.$$
Passing to the infimum over $u$, we have $\inf_{\mc{P}_{(\rho)}} I^m \leq \inf_{\mc{P}_{\mu}} J_{\mu} - \frac{\mu}{2} m$. Passing to the infimum over $\mu$, we have $\inf_{\mc{P}_{(\rho)}} I^m \leq \inf_{
\rho < \mu < \overline{\mu}_0} \left(\inf_{\mc{P}_{\mu}} J_{\mu}-\frac{\mu}{2}m\right)$.

Let now $(\nu, u) \in\mc{P}_{(\rho)
}$, which implies 
$\nu \in (
\rho, \overline{\mu}_0)$. As a consequence
$$I^m(\nu, u) =J_{\nu}(u) - \frac{\nu}{2} m \geq \inf_{\mc{P}_{\nu}} J_{\nu} - \frac{\nu}{2} m \geq \inf_{\rho < \mu < \overline{\mu}_0} \left(\inf_{\mc{P}_{\mu}} J_{\mu}-\frac{\mu}{2}m\right).$$
Passing to the infimum over $(\nu,u)$ we obtain $\inf_{\mc{P}_{(\rho)}} I^m \geq \inf_{
\rho < \mu < \overline{\mu}_0} \left(\inf_{\mc{P}_{\mu}} J_{\mu}-\frac{\mu}{2}m\right)$,
whence the claim.
\QED




\subsection{$L^2$-minima, Poho\v{z}aev identity, and critical points}
\label{subsub_l2_poh_crit}

\begin{Proposition}\label{pr:PohoMin}
	Assume \hyperref[g0]{\rm{(g0)}}--\hyperref[g2]{\rm{(g2)}}. 
Let $u \in H^s(\R^N)$.
\begin{itemize}
\item[(i)] If $u \in \mc{S}_{m}
$ is an $L^2$-minimum for $K$, i.e., $K(u)= \inf_{\mc{S}_m} K 
$, 
then $u$ is a 
solution of \eqref{eq_general_mu_fixed} for some Lagrange multiplier $\mu \in \R$, and $(\mu,u)$ satisfies the Poho\v{z}aev identity. 

\item[(ii)] If $u \in \mc{D}_{m}
$ is such that $K(u)= \inf_{\mc{D}_m} K $ and $|u|_2^2<m$, then $u$ is a 
solution of \eqref{eq_general_mu_fixed} with Lagrange multiplier $\mu=0$, and $(\mu,u)$ satisfies the Poho\v{z}aev identity.

\item[(iii)] If $(\mu,u) \in \mc{P}$
and $K(u) \le 0$, then $\mu>0$.
\item[(iv)]
If $u \in \mc{D}_{m}
$ is such that $K(u)= \inf_{\mc{D}_m} K <0$, then $u \in \mc{S}_m$; in particular, $K(u)= \inf_{\mc{S}_m} K$. 
\end{itemize}
\end{Proposition}

\claim Proof.
(i) First, we show that $u$ is a 
solution. Indeed, let $\varphi \in H^s
(\R^N) \cap L^1(\R^N)$; being $u \neq 
 0$, for $t \in \R$ with $|t|$ small we can assume $u + t \varphi \neq 
 0$. Thus, for each $t$ we choose $\theta(t) > 0$ such that $\abs{ (u+t \varphi)(\cdot/\theta(t))}_2^2=m$, that is,
$\theta(t):= \left( m\abs{u+t \varphi}_2^{-2}\right)^{1/N}.$ 
Let $\gamma(t):= K\big((u+t \varphi)(\cdot/\theta(t))\big)$; due to the differentiability of $\theta$ and the arguments of the proof of Lemma \ref{lem_diff_along_L1}, we have that $\gamma$ is differentiable.
Being $t=0$ a point of minimum for $\gamma$, 
we have $\gamma'(0)=0$, that is (notice $\theta(0)=1$),
$$\frac{N-2s}{2}\theta'(0)\abs{(-\Delta)^{s/2}u}_2^2 + \intRN(-\Delta)^{s/2}u (-\Delta)^{s/2}\varphi \, \dx - N\theta'(0)\intRN G(u) \, \dx - \intRN g(u)\varphi \, \dx = 0$$
and hence (being $\theta'(t) = - \frac{2}{N}m^{1/N} \abs{u+t\varphi}_2^{-2/N-2} \intRN (u+t \varphi) \varphi \, \dx$)
$$
\intRN(-\Delta)^{s/2}u (-\Delta)^{s/2}\varphi \, \dx - \frac{2}{m} \left(\frac{1}{2^*_s}\abs{(-\Delta)^{s/2}u}_2^2 - \intRN G(u) \, \dx\right) \intRN u \varphi \, \dx - \intRN g(u)\varphi \, \dx = 0.
$$
Set 
\begin{equation}\label{eq_def_lagrang_mult}
\mu:=-\frac{2}{m} \left(\frac{1}{2^*_s}\abs{(-\Delta)^{s/2}u}_2^2 - \intRN G(u)\right) \in \R,
\end{equation}
 we have the claim. 
We see that the Poho\v{z}aev identity directly comes from the definition of $\mu$ in \eqref{eq_def_lagrang_mult}.
%

(ii) The statement is obvious for $u\equiv 0$. Arguing as in (i), let $\varphi \in H^s
(\R^N) \cap L^1(\R^N)$; being $u \neq 
 0$, for $t \in \R$ with $|t|$ small we can assume $u + t \varphi \neq 
 0$ and $|u+t \varphi|_2^2 <m$. By the minimality of $t=0$ for $t \mapsto K(u+t\varphi)$ we obtain
$$ \intRN(-\Delta)^{s/2}u (-\Delta)^{s/2}\varphi \, \dx - \intRN g(u)\varphi \, \dx = 0$$
which implies that $u$ is a solution with zero Lagrange multiplier. 
To show the Poho\v{z}aev identity, consider the map $t \mapsto K(u(\cdot/t))$, where $|u(\cdot/t)|_2^2<m$ for $t$ close to $1$. Being $t=1$ a point of minimum, we have
$$\frac{N-2s}{2} |(-\Delta)^{s/2} u|_2^2 - N \int_{\R^N} G(u) \, \dx = 0$$
which is the claim.

The proof of (iii) is straightforward. Point (iv) follows by (ii) and (iii).
\QED

\bigskip

We prove now Proposition \ref{pr-Poho_sol}. We underline that, unlike \cite{BKS17}, we employ a direct proof without relying on a deformation lemma.


\medskip

\claim Proof of Proposition \ref{pr-Poho_sol}.
Let $r,t \in \R$ and $\varphi \in H^s(\R^N) \cap L^1(\R^N)$, and define $\sigma = \sigma(r,t) > 0$ such that $P\bigl(\mu+r,u(\cdot/\sigma) + t \varphi(\cdot/\sigma)\bigr) = 0$, that is,
\[
\sigma(r,t) := \left(\frac{1}{2^*_s} |(-\Delta)^{s/2}(u+t\varphi)|_2^2 \left(\intRN G(u + t\varphi) - \frac{\mu+r}{2} (u + t\varphi)^2 \, \dx\right)^{-1}\right)^{\frac{1}{2s}}.
\]
Indeed, since $P(\mu,u) = 0$, we know that $\intRN G(u) - \frac{\mu}{2} u^2 \, \dx>0$, thus for $|r|, |t|$ small the above expression is well defined.
Moreover we can assume $|r|$ small such that $\mu+r > \rho$, thus $(\mu+r,u(\cdot/\sigma) + t \varphi(\cdot/\sigma)) \in \mc{P}_{(\rho)}$.

Being $P(\mu,u) = 0$, we have $\sigma(0,0) = 1$. Next, we compute
\begin{align*}
\MoveEqLeft I^m\bigl(\mu+r,u(\cdot/\sigma) + t \varphi(\cdot/\sigma)\bigr) \\
&= \frac{\sigma^{N-2s}}{2} |(-\Delta)^{s/2}(u + t\varphi)|_2^2 + \frac{\mu + r}{2} \big(\sigma^N |u + t\varphi|_2^2 - m\big) - \sigma^N \intRN G(u + t\varphi) \, \dx.
\end{align*}
By Lemma \ref{lem_diff_along_L1}, the map $t \mapsto \intRN G(u + t\varphi) \, \dx$, and so $\sigma$ as well, is differentiable at $0$. Since $(\mu,u)$ is a minimizer, the gradient of $I^m\bigl(\mu+r,u(\cdot/\sigma) + t \varphi(\cdot/\sigma)\bigr)$ at $(r,t) = (0,0)$ is zero. At the same time, recalling that $\sigma(0,0) = 1$,
\begin{align*}
\partial_r I^m\bigl(\mu+r,u(\cdot/\sigma)\bigr)|_{(t,r)=(0,0)} = \frac{N-2s}{2} \partial_r \sigma(0,0) P(\mu,u) + \frac12 (|u|_2^2 - m) = \frac12 (|u|_2^2 - m)
\end{align*}
and
\begin{align*}
\MoveEqLeft \partial_t I^m\bigl(\mu,u(\cdot/\sigma) + t \varphi(\cdot/\sigma)\bigr)|_{(t,r)=(0,0)} = \\
&= \frac{N-2s}{2} \partial_t \sigma(0,0) P(\mu,u) + \int_{\R^N} 
(-\Delta)^{s/2} u (-\Delta)^{s/2} \varphi + \mu u \varphi - g(u) \varphi 
 \, \dx\\
& = \int_{\R^N} 
(-\Delta)^{s/2} u (-\Delta)^{s/2} \varphi + \mu u \varphi - g(u) \varphi 
\, \dx.
\end{align*}
Thus, we obtain the claim. We conclude by applying Proposition \ref{prop_uguag_not_sym}.
\QED

\bigskip

\claim Proof of Proposition \ref{prop_intr_L2min_PN}. 
It follows by Proposition \ref{pr:PohoMin} and Proposition \ref{prop_uguag_not_sym}.
%
\QED

\bigskip

\claim Proof of Proposition \ref{prop_JJLU}.

(i) 
From Proposition \ref{prop_intr_L2min_PN} we have that $I^m(\mu,u) = \inf_{\mc{P}} I^m$, hence, by Remark \ref{rem_2min}, $J_{\mu}(u) 
=\inf_{\mc{P}_{\mu}} J_{\mu}$, and \eqref{eq_relat_kep_a} follows from Proposition \ref{prop_uguag_not_sym}.

(ii) By \eqref{eq_relat_kep_a} and Proposition \ref{prop_uguag_not_sym} we have
$$J_{\mu}(u) 
= \inf_{\mc{P}_{\mu}} J_{\mu} = \inf_{\mc{S}_m} K + \frac{\mu}{2} m = \inf_{\mc{P}} I^m + \frac{\mu}{2}m,$$
that is, $I^m(\mu, u) = \inf_{\mc{P}} I^m$. From Proposition \ref{pr-Poho_sol} 
we obtain that $u \in \mc{S}_m$ and we conclude.
\QED

\subsection{Nehari identity}

\begin{Proposition}[Nehari identity]\label{p-Nehari}
Assume \hyperref[g0]{\rm{(g0)}}--\hyperref[g2]{\rm{(g2)}} and let $\mu \in \R$.
If $u \in H^s(\R^N)$ is a 
solution to \eqref{eq_general_mu_fixed} such that $\intRN g_-(u) u \, \dx < \infty$, 
then $u$ satisfies the Nehari identity
\[
|(-\Delta)^{s/2} u|_2^2 + \mu |u|_2^2 = \intRN g(u) u \, \dx.
\]
\end{Proposition}
\claim Proof.
With the notation of Lemma \ref{lem_converg_Hs}, for every $R>0$ 
there holds
\begin{equation}\label{e-LMAO_2}
\intRN (-\Delta)^{s/2}u \ (-\Delta)^{s/2}(\varphi_R u) \, \dx + \mu \intRN u^2 \varphi_R \, \dx = \intRN g(u) \varphi_R u \, \dx.
\end{equation}
Since $|\varphi_R| \le 1$ and $\lim_{R\to +\infty} \varphi_R = 1$ a.e. in $\R^N$, from the dominated convergence theorem we obtain
\[
\lim_{R\to +\infty} \intRN u^2 \varphi_R \, \dx = \intRN u^2 \, \dx \quad \text{and} \quad \lim_{R\to +\infty} \intRN g(u) \varphi_R u \, \dx = \intRN g(u) u \, \dx;
\]
using also Lemma \ref{lem_converg_Hs}, we can pass to the limit in \eqref{e-LMAO_2} and conclude. 
\QED

%

\section{The superlinear (perturbed) problem}\label{s-pert}


\indent

In this section, 
where we always assume \hyperref[g0]{\rm{(g0)}}--\hyperref[g2]{\rm{(g2)}}, we study the following perturbed problem 
\begin{equation}\label{e-pert}
\begin{cases}
\displaystyle (-\Delta)^s u + \mu u = g_\eps(u),\\
\displaystyle \intRN u^2 \, \dx = m,\\
(\mu,u) \in (0,+\infty) \times H^s(\R^N),
\end{cases}
\end{equation}
where, for $\eps \in (0,1]$, $g_\eps$ is defined as follows:
let $\varphi_{\eps} \in \mc{C}(\R, [0,1])$ be given by
\begin{equation}\label{eq_def_varphieps}
\varphi_{\eps}(t) :=
\begin{cases}
|t|/\eps & \quad \hbox{for $|t|< \eps$}, \\
1 & \quad \hbox{for $|t|\geq \eps$},
\end{cases}
\end{equation}
and set
\begin{equation*}
g_{-}^{\eps}(t):= \varphi_{\eps}(t) g_-(t), \quad G_-^{\eps}(t):=\int_0^t g_-^{\eps}(\tau) \, \mathrm{d}\tau.
\end{equation*}
Then, we define
\begin{equation}\label{eq_def_geps}
g_{\eps}(t):= g_+(t) - g_{-}^{\eps}(t), \quad G_{\eps}(t):=\int_0^t g_{\eps}(\tau) \, \mathrm{d}\tau = G_{+}
(t) - G_-^{\eps}(t).
\end{equation}
Notice that $\varphi_{\eps} \to 1$ pointwise, thus
$$ g_{\eps}(t)\to g(t), \quad G_{\eps}(t) \to G(t) \quad \hbox{for every $t \in \R$}$$
as $\eps \to 0^+$. 
Hence, \eqref{e-pert} can be considered as an approximating family of problems; the advantage is given by the modification in zero, which suppresses the singularity in the sublinear case: 
\begin{equation}\label{eq_extra_cond_eps}
\lim_{t \to 0} \frac{g_{\eps}(t)}{t} =0.
\end{equation}
%
%
%
%
%
%
We will make large use of the following monotonicity property: if $0 < \eps' < \eps \le 1$, then 
$$G_-^\eps \le G_-^{\eps'} \leq G_-,$$
which implies, moreover,
$$G_\eps \ge G_{\eps'} \ge G. 
$$

\medskip

The purpose of this section is 
to obtain the existence of a Poho\v{z}aev minimum for the problem \eqref{e-pert}, furnishing, in addition, a slight generalization of \cite{CGT21N} by exploiting the refined Theorem \ref{t-bks}; such a generalization holds for more general nonlinearities not related to the perturbation process, see Remark \ref{t-g*} below.



\subsection{$\eps$ fixed: some refinement for the superlinear problem}
\label{subsec_minim_pohoz}

\indent

We start by considering $\eps \in (0,1]$ fixed. Assuming $\mu>0$, to avoid technical issues related to the boundary of $(0,+\infty)$, we exploit the identification
\begin{equation}\label{eq_identif_lambda_mu}
\mu \equiv e^{\lambda}.
\end{equation}
A similar approach that does not make use of \eqref{eq_identif_lambda_mu} is developed in \cite[Section 4.2]{Gal23}.

\begin{Remark}
The condition $\mu>0$ is a classical requirement for nonlinearities with superlinear growth at the origin -- cf. \eqref{eq_extra_cond_eps} -- because it leads to (see \cite{BeLi83I,ChWa13})
\[
\lim_{t \to 0} \frac{g_\eps(t) - \mu t}{t} < 0.
\]
\end{Remark}
Making use of \eqref{eq_identif_lambda_mu}, we can thus reformulate the search for 
solutions to \eqref{e-pert} as 
critical points of the Lagrangian functional $I^m_{\eps} \colon \R \times H^s(\R^N) \to \R$ defined as
\begin{equation*}
I^m_{\eps}(\lambda,u) := \frac12 |(-\Delta)^{s/2} u|_2^2 + \frac{e^\lambda}{2} (|u|_2^2 - m) - \intRN G_\eps(u) \, \dx. 
\end{equation*}
Let us also define $P_\varepsilon \colon \R \times H^s(\R^N) \to \R$ as
\begin{equation*}
P_{\eps}(\lambda,u) := |(-\Delta)^{s/2} u|_2^2 + 2^*_s \left(\frac{e^\lambda}{2} |u|_2^2 - \intRN G_\eps(u) \, \dx\right),
\end{equation*}
\[
\mc{P}_{\eps,+} := \big\{(\lambda,u) \in \R \times H^s(\R^N) \mid P_\eps(\lambda,u) = 0, \, u \ne 0\big\}, \quad \mc{P}_{\eps,+}^{\textup{rad}}:= \mc{P}_{\eps,+} \cap \bigl(\R \times H^s_{\textup{rad}}(\R^N)\bigr),
\]
and
\begin{equation*}
d_{\eps,+}^m := \inf_{\mc{P}_{\eps,+}^{\textup{rad}}} I^m_{\eps}.
\end{equation*}
With an abuse of notation, we redefine consequently on $\R \times H^s(\R^N)$ also
\begin{equation*}
I^m(\lambda, u), \; P(\lambda, u)
\end{equation*}
defined in \eqref{eq_def_Im} and \eqref{eq_def_Poh_pos}.
We also 
set
\begin{equation}\label{eq_def_P+}
\mc{P}_{+} := \big\{(\lambda,u) \in \R \times H^s(\R^N) \mid P(\lambda,u) = 0, \, u \ne 0\big\}, \quad \mc{P}_{+}^{\textup{rad}}:= \mc{P}_{+} \cap \left(\R \times H^s_{\textup{rad}}(\R^N)\right).
\end{equation}
Notice that, with the identification $\mu \equiv e^\lambda > 0$, we have
\begin{equation}\label{eq_def_em+}
d_+^m:= \inf_{(\lambda,u) \in \mc{P}_{+}^\textup{rad}} I^m(\lambda,u)= \inf_{(\mu,u) \in \mc{P}_{(0)}^\textup{rad}} I^m(\mu,u) .
\end{equation}
Moreover, we introduce
\begin{equation*}
K_{\eps}(u) := \frac12 |(-\Delta)^{s/2} u|_2^2 - \intRN G_\eps(u) \, \dx, \quad \kappa_{\eps}^m := \inf_{\mc{S}_m^{\textup{rad}}} K_{\eps}, \quad \ell_{\eps}^m := \inf_{\mc{D}_m^{\textup{rad}}} K_{\eps}.
\end{equation*}
The following is the main theorem of this section.

\begin{Theorem}\label{t-pert}
Assume that \hyperref[g0]{\rm{(g0)}}--\hyperref[g4]{\rm{(g4)}} hold, then there exists $m_\eps \geq 0$ such that for every $m > m_\eps$ there exists a solution $(\lambda
_\eps,u_\eps) \in \mc{P}_{\eps,+}^{\textup{rad}}$ to \eqref{e-pert} verifying $I^m_{\eps}(\lambda_{\eps}, u_{\eps})=d_{\eps,+}^m<0$. 

Moreover, if $g|_{(-\infty,0)} = 0$ or $g$ is odd, then for every $m > m_\eps$ there exists a Schwarz-symmetric minimizer of $I^m_{\eps}$ over $\mc{P}_{\eps,+}$;
in particular,
\[
\inf_{\mc{P}_{\eps,+}
} I^m_{\eps} = \inf_{\mc{P}_{\eps,+}^{\textup{rad}} 
} I^m_{\eps}.
\]
\end{Theorem}

\begin{Remark}\label{t-g*}
We highlight that the only additional property of $g_\eps$ with respect to \hyperref[g0]{\rm{(g0)}}--\hyperref[g4]{\rm{(g4)}} that is used in Theorem \ref{t-pert} is the superlinearity at zero of $g_{\eps}$, i.e., \eqref{eq_extra_cond_eps} (or equivalently, from \hyperref[g1]{\rm{(g1)}}, $\lim_{t \to 0} g_-^\eps(t) / t = 0$). Consequently, Theorem \ref{t-pert} and all the other properties presented in this section 
still hold true with $g$ instead of $g_\eps$ under the additional assumption
\[
\lim_{t \to 0}\frac{g_-(t)}{t} = 0.
\]
This furnishes a small generalization of the results in \cite{CGT21N}.
\end{Remark}


\begin{Remark}\label{rem_multipl}
We underline that the solution found in Theorem \ref{t-pert} 
has also a minimax characterization.
Introduced 
\begin{equation*}
	\Omega_{\eps} := \left\{(\lambda,u) \in \R \times H^s(\R^N) \mid P_\eps(\lambda,u) > 0 \right\} \cup \big(\R \times \{0\}\big),
\end{equation*}
we have $\partial \Omega_{\eps} = \mc{P}_{\eps,+}$ (see \cite[Lemma 2.3.1]{Gal23}),
and this is treated as a ``mountain'' in the minimax approach: 
set
\begin{equation*}
\Gamma^m_{\eps} := \left\{\xi \in \mc{C}\bigl([0,1],\R \times H^s_{\textup{rad}}(\R^N)\bigr) \;\middle|\; \xi(0) \in \R \times \{0\}, \; \xi(1) \not\in \Omega_{\eps}, \; \max_{i=0,1} I^m_{\eps}\bigl(\xi(i)\bigr) \le d_{\eps,+}^m - 1\right\},
\end{equation*}
one can prove that \cite[Section 6]{CGT21N}
$$d_{\eps,+}^m =\inf_{\xi \in \Gamma^m_{\eps}} \max_{t \in [0,1]} I^m_{\eps}\bigl(\xi(t)\bigr).$$
Moreover, by arguing 
as in \cite{CGT21N}, 
we also obtain (write $g_0 := g$ and $G_0 := G$ for convenience)
\begin{itemize}
	\item if $g$ is odd, then for each $\eps \in (0,1]$ and $k \in \N$ there exists $m_{\eps,k} \geq 0$ such that, for $m >m_{\eps,k}$, there exist $k$ (pairs 
	of) solutions to \eqref{e-pert};
	\item if $\displaystyle \lim_{t \to 0} \frac{G(t)}{|t|^{\bar{p}+1}}=+\infty$, then for each $\eps \in [0,1]$, we have $m_{\eps}=0$ (and $m_{\eps,k}=0$ if $g$ is odd) because $G_\eps \ge G$.
\end{itemize}
\end{Remark}

In \cite{CGT21N}, the existence of a normalized solution is pursued by a Mountain Pass approach, involving a suitable augmented functional and a deformation lemma; the found solution is shown a posteriori to be a Poho\v{z}aev minimum. 
Here we provide a different approach, by minimizing directly the functional on the Poho\v{z}aev product set: this proof is not trivial, since a precise sequence of minimizers has to be selected; on the other hand, it requires fewer tools. 

To prove Theorem \ref{t-pert}, we need 
preliminary estimates on $d_{\eps,+}^m$. This will be given by the interplay of the value of the Poho\v{z}aev minimum on the product space and the one with fixed $\mu$.


\subsubsection{The Poho\v{z}aev asymptotic geometry} 

\indent

To estimate $d_{\eps,+}^m$, we introduce now the Poho\v{z}aev geometry 
related to a fixed frequency $\lambda$, in the spirit of Section \ref{s-BKS}. 
 For $\lambda \in \R$ we define 
\begin{equation}\label{e-mpGamma}
\begin{split}
J(\lambda,u) &: = \frac12 |(-\Delta)^{s/2} u|_2^2 + \frac{e^\lambda}{2} |u|_2^2 - \intRN G(u) \, \dx, 
\\
J_\varepsilon(\lambda,u) &: = \frac12 |(-\Delta)^{s/2} u|_2^2 + \frac{e^\lambda}{2} |u|_2^2 - \intRN G_\eps(u) \, \dx, 
\\
\mc{P}(\lambda) &:= \left\{u \in H^s(\R^N) \mid (\lambda,u) \in \mc{P}_+ \right\}, 
\quad \mc{P}^{\textup{rad}}(\lambda) := \mc{P}(\lambda) \cap H^s_{\textup{rad}}(\R^N), \\
 \mc{P}_{\eps}(\lambda) & := \left\{u \in H^s(\R^N) \mid (\lambda,u) \in \mc{P}_{\eps,+} \right\}, \quad \mc{P}_{\eps}^{\textup{rad}}(\lambda) := \mc{P}_{\eps}(\lambda) \cap H^s_{\textup{rad}}(\R^N), \\
a(\lambda) &:= \inf_{\mc{P}^{\textup{rad}}(\lambda)} J(\lambda, \cdot), 
\quad a_\eps(\lambda) := \inf_{\mc{P}_{\eps}^{\textup{rad}}(\lambda)} J_{\eps}(\lambda, \cdot) .
\end{split}
\end{equation}
Notice that $J(\lambda, u)=J_{e^{\lambda}}(u)$ and $\mc{P}(\lambda)=\mc{P}_{e^{\lambda}}$ 
with the notations of Section \ref{s-BKS}. Moreover, we make again an abuse of notation by writing $a(\lambda) \equiv a(e^{\lambda})$ -- see \eqref{eq_def_amu}.
We observe that 
$$J_{\eps} \le J.$$
To compare $a(\lambda)$ with $a_{\eps}(\lambda)$, we make use of \eqref{eq_min_poh_less}: as a matter of fact, by $P_{\eps}\leq P$ we have $\{P=0\} \subset \{P_{\eps}\leq 0\}$ and thus 
\begin{equation}\label{eq_comp_aeps_a0}
a_\eps(\lambda) \le a(\lambda).
\end{equation}
%
We highlight that, arguing as in \cite{BKS17}, we could give a Mountain Pass characterization of $a(\lambda)$; 
nonetheless, we will not make use of this characterization for the limit problem.
Moreover, we introduce
\begin{equation}\label{e-mu0}
 \begin{split}
 \overline{\mu}_\eps := & \, \sup_{t\ne0} \tfrac{G_\eps(t)}{t^2/2}\\
 \geq & \,\sup_{t\ne0} \tfrac{G(t)}{t^2/2} = \overline{\mu}_0 \in (0, +\infty] \\
 \end{split}
\qquad
 \begin{split}
 \overline{\lambda}_\eps := & \, \log(\overline{\mu}_\eps)\\
 \geq & \, \log(\overline{\mu}_0)=:\overline{\lambda}_0 \in (-\infty, +\infty]
 \end{split}
\end{equation}
%
%
%
%
with the convention that $\log(+\infty) = +\infty$. Observe that $\overline{\mu}_0 > 0$ follows from \hyperref[g4]{\rm{(g4)}}. We highlight that, in what follows, we will sometimes need to distinguish the cases $\overline{\lambda}_\eps \in \R$ and $\overline{\lambda}_\eps = +\infty$.

We define, for $m>0$, the values
\begin{equation}\label{e-infima}
b_{\eps}^m := \inf_{\lambda < \overline{\lambda}_\eps} \left(a_\eps(\lambda) - \frac{e^\lambda}{2} m\right) \quad \text{and} \quad b^m := \inf_{\lambda < \overline{\lambda}_0} \left(a(\lambda) - \frac{e^\lambda}{2} m\right).
\end{equation}
Moreover, by Lemma \ref{l-lambda-gen} the following thresholds are well defined
\begin{equation}\label{eq_meps_m0}
m_\eps := \inf_{\lambda < \overline{\lambda}_\eps} \frac{a_\eps(\lambda)}{e^\lambda/2} \quad \text{and} \quad m_0 := \inf_{\lambda < \overline{\lambda}_0} \frac{a(\lambda)}{e^\lambda/2}.
\end{equation}
Observe that, from Lemma \ref{l-lambda-gen} and since $\overline{\lambda}_\eps \ge \overline{\lambda}_0$ and $a 
\geq a_\eps \geq 0$ -- cf. Proposition \ref{p-a_monot}, 
we have 
\begin{equation}\label{eq_compar_meps_m0}
b^m \geq b_{\eps}^m \quad \text{and} \quad m_0 \ge m_\eps \ge 0.
\end{equation}


%


As a straightforward consequence of Proposition \ref{prop_intr_legendre}, Theorem \ref{t-bks}, and the definitions of $\overline{\lambda}_{\eps}$ and $m_{\eps}$, we obtain the following result. 
\begin{Proposition}\label{p-sol2free}
Assume \hyperref[g0]{\rm{(g0)}}--\hyperref[g2]{\rm{(g2)}} and \hyperref[g4]{\rm{(g4)}}. Then 
\begin{itemize}
\item For every $\lambda < \overline{\lambda}_{\eps}$
there exists a solution $u_{\lambda,\eps} \in \mc{P}_{\eps}^{\textup{rad}}(\lambda) 
$ to
$$(-\Delta)^s u + e^{\lambda} u = g_\eps(u) \quad \hbox{in $\R^N$}.$$
Additionally,
\[
J_{\eps}(\lambda,u_{\lambda,\eps}) = \inf_{\mc{P}_{\eps}^{\textup{rad}}(\lambda)}J_{\eps}(\lambda,u)
= a_\eps(\lambda)>0.
\]
As a consequence, $a(\lambda)>0$.
\item $d_{\eps,+}^m = b_{\eps}^m$; if $m>m_{\eps}$, then $d_{\eps,+}^m<0$. 
\item If $\overline{\lambda}_\eps \in \R$, then $\min\{P_{\eps}(\lambda,u),J_{\eps}(\lambda,u)\} > 0$ for all $\lambda \ge \overline{\lambda}_\eps$ and $u \in H^s(\R^N) \setminus \{0\}$.
\end{itemize}
\end{Proposition}

We focus now on properties strictly related to the $L^2$-subcritical setting \hyperref[g3]{\rm{(g3)}}.
To start, we observe that for every $\delta>0$ there exist $C_\delta,c_\delta>0$ such that for every $t \in \R$
\begin{equation}\label{eq_stima_Geps_1}
G_\eps(t) \le G_+(t) \le \frac{\delta}{2} 
t^2 + \frac{C_\delta}{\bar{p}+1}
 |t|^{\bar{p}+1}
\end{equation}
and
\begin{equation}\label{e-blue}
G_\eps(t) \le G_+(t) \le \frac{ c_\delta}{2} 
t^2 + \frac{\delta}{\bar{p}+1}
 |t|^{\bar{p}+1}.
\end{equation}

\begin{Lemma}[Asymptotic geometry]
\label{l-bddB}
Assume \hyperref[g0]{\rm{(g0)}}--\hyperref[g4]{\rm{(g4)}}.
The following statements hold.
\begin{itemize}
\item 
$\lim_{\lambda \to \overline{\lambda}_\eps^-} a_\eps(\lambda) = +\infty$. 
If $\overline{\lambda}_\eps = +\infty$, then $\lim_{\lambda \to +\infty} \frac{a_\eps(\lambda)}{e^\lambda} = +\infty$.
\item 
$d_{\eps,+}^m 
>-\infty$. 
\end{itemize}
\end{Lemma}


\claim Proof.
Recalled \eqref{eq_stima_Geps_1}-\eqref{e-blue}, concerning the first point we have that
the first limit can be obtained with the same reasoning as in \cite[Proposition 3.3]{CGT21N} in 
light of the monotonicity in Proposition \ref{p-a_monot}, while for the second one we can argue as in \cite[Proposition 3.4]{CGT21N} once the following property is proved: 
\begin{equation}\label{e-delta}
a_\varepsilon(\lambda) \ge (e^\lambda - c_\delta) \inf_{\mc{P}^*(\delta)} J_\delta^* \; \; \text{ if } e^\lambda > c_\delta \quad \text{and} \quad \lim_{\delta \to 0^+} \inf_{\mc{P}^*(\delta)} J_\delta^* = +\infty,
\end{equation}
where
\begin{align*}
J_\delta^*(u) & := \frac12 |(-\Delta)^{s/2}u|_2^2 + \frac12 |u|_2^2 - \frac{\delta}{\bar{p}+1} |u|_{\bar{p}+1}^{\bar{p}+1}, \\
\mc{P}^*_{\delta} 
 & := \left\{ u \in H^s_\textup{rad}(\R^N) \mid P_\delta^*(u) = 0 \right\} \setminus \{0\}, \\
P_\delta^*(u) & := |(-\Delta)^{s/2}u|_2^2 + 2_s^* \left( \frac12 |u|_2^2 - \frac{\delta}{\bar{p}+1} |u|_{\bar{p}+1}^{\bar{p}+1} \right)
\end{align*}
(i.e., $P_\delta^*$ and $\mc{P}^*_{\delta} 
$ are, respectively, the Poho\v{z}aev functional and manifold related to $J_\delta^*$). 
For every $\theta > 0$ and $u \in \mc{P}_{\eps}^{\textup{rad}}(\lambda)$, we define $u_\theta := \theta^{-N/(4s)} u(\theta^{-1/(2s)} \cdot)$. Then, using \eqref{e-blue} and taking $\theta = e^\lambda - c_\delta$ (which is positive if $\lambda > \log (c_\delta)$), it is straightforward to verify that
\[
J_{\eps}(\lambda,u) \ge \theta J_\delta^*(u_{\theta}) \quad \text{and} \quad P_\delta^*(u_{\theta}) \le 0, 
\]
therefore the first property in \eqref{e-delta} holds thanks to \eqref{eq_min_poh_less}. 
As for the second one, by scaling, one readily checks that $P_\delta^*(\delta^{-1/(\bar{p}-1)} \cdot) = \delta^{-2/(\bar{p}-1)} P_1^*$, which implies that $\mc{P}^*_{\delta} = \delta^{-1/(\bar{p}-1)} \mc{P}^*_1$, and $J_\delta^*(\delta^{-1/(\bar{p}-1)} \cdot) = \delta^{-2/(\bar{p}-1)} J_1^*$. 
From this, one obtains $\inf_{\mc{P}^*_{\delta}} J_\delta^* = \delta^{-2/(\bar{p}-1)} \inf_{\mc{P}^*_1} J_1^* \to +\infty$ as $\delta \to 0^+$.

Finally, recalled that $d_{\eps,+}^m = b_{\eps}^m $ by Proposition \ref{p-sol2free}, the estimate $b_{\eps}^m > -\infty$ is obtained as in \cite[Proposition 3.7]{CGT21N}. 
\QED
\subsubsection{Direct minimization on Poho\v{z}aev product}


\indent

We are now ready to prove the existence of a Poho\v{z}aev minimum as stated in Theorem \ref{t-pert}. 

\smallskip


\claim Proof of Theorem \ref{t-pert}. 
By Proposition \ref{p-sol2free} and Lemma \ref{l-bddB} 
we know that $d_{\eps,+}^m \in (-\infty, 0)$. Let $(\lambda_n, u_n)_n \subset \mc{P}_{\eps,+}^{\textup{rad}}$ 
be a minimizing sequence, i.e., $I^m_{\eps}(\lambda_n, u_n) \to d_{\eps,+}^m$ and $P_{\eps}(\lambda_n, u_n)=0$. 
Notice that, being $u_n \in \mc{P}_{\eps}^{\textup{rad}}(\lambda_n)$, 
by Lemma \ref{l-lambda-gen} we have $\lambda_n < \overline{\lambda}_{\eps}$ for each $n$.

\smallskip
\textbf{Step 1.} \emph{Choice of a proper minimizing sequence}.

We build now a minimizing sequence with additional information on the energy levels and on the differential (with respect to $u$) of $I_\varepsilon^m$.
%
For each $n$, by 
Proposition \ref{p-sol2free}
there exists $v_n \in \mc{P}_{\eps}^{\textup{rad}}(\lambda_n)$ 
such that 
\begin{equation}\label{eq_dim_extra_prop1}
I^m_{\eps}(\lambda_n, v_n) = a_\eps(\lambda_n) - \frac{e^{\lambda_n}}{2} m
\end{equation}
and
\begin{equation}\label{eq_dim_extra_prop2}
\partial_u I^m_{\eps}(\lambda_n, v_n)=0.
\end{equation}
From \eqref{eq_dim_extra_prop1} we obtain
$$d_{\eps,+}^m \leq I^m_{\eps}(\lambda_n, v_n) = a_\eps(\lambda_n) - \frac{e^{\lambda_n}}{2} m \leq J_{\eps}(\lambda_n, u_n) - \frac{e^{\lambda_n}}{2}m =I^m_{\eps}(\lambda_n,u_n) =d_{\eps,+}^m + o(1).$$
As a consequence, also $(\lambda_n, v_n)_n$ is a minimizing sequence. 

\smallskip
\textbf{Step 2.} \emph{$(\lambda_n)_n$ is bounded below}.

We have
\begin{equation*}
\frac{m}{2} e^{\lambda_n} = \frac{1}{2^*_s} P_\eps(\lambda_n,v_n) - I^m_{\eps}(\lambda_n,v_n) + \frac{s}{N} |(-\Delta)^{s/2} v_n|_2^2 \ge
- I^m_{\eps}(\lambda_n,v_n),
\end{equation*}
whence
\begin{equation*}
\liminf_n \frac{m}{2} e^{\lambda_n} \ge - d_{\eps,+}^m > 0.
\end{equation*}

\smallskip
\textbf{Step 3.} \emph{$(\lambda_n)_n$ is bounded above}.

In this Step we use the $L^2$-subcritical nature of the problem.
If $\overline{\lambda}_{\eps}<+\infty$, then we are done. Otherwise, if by contradiction $\lambda_n \to +\infty$, then by Lemma \ref{l-bddB} and \eqref{eq_dim_extra_prop1} we have
$$d_{\eps,+}^m + o(1)= I^m_{\eps}(\lambda_n, v_n) = e^{\lambda_n} \left( \frac{a_\eps(\lambda_n)}{e^{\lambda_n}} -\frac{1}{2} m\right) \to +\infty,$$
which is impossible. 

\smallskip
\textbf{Step 4.} \emph{$(|(-\Delta)^{s/2} v_n|_2)_n$ is bounded}.


Indeed,
\begin{equation*}
\frac{s}{N} |(-\Delta)^{s/2} v_n|_2^2 = I^m_{\eps}(\lambda_n,v_n) - \frac{1}{2^*_s} P_\eps(\lambda_n,v_n) + \frac{m}{2} e^{\lambda_n}= d_{\eps,+}^m + \frac{m}{2} e^{\lambda_n} + o(1)
\end{equation*}
and thus $|(-\Delta)^{s/2} v_n|_2^2$ is bounded thanks to Step 3. 

\smallskip
\textbf{Step 5.} \emph{$(|v_n|_2)_n$ is bounded}.

By \eqref{eq_stima_Geps_1} we have
\begin{align*}
d_{\eps,+}^m + \frac{m}{2} e^{\lambda_n} + o(1)
	& = I^m_{\eps}(\lambda_n,v_n) - \frac{1}{2} P_\eps(\lambda_n,v_n) + \frac{m}{2} e^{\lambda_n} \\
	&= -\frac{e^{\lambda_n}}{2} \left(\frac{2^*_s}{2} - 1\right) |v_n|_2^2 + \left(\frac{2^*_s}{2} - 1\right) \intRN G_{\eps}(v_n) \, \dx\\
	&  \le \left(\frac{2^*_s}{2} - 1\right) \left( -\frac{1}{2}\left(e^{\lambda_n} -\delta\right) |v_n|_2^2 + \frac{C_{\delta}}{\bar{p}+1} |v_n|_{2^*_s}^{2^*_s} \right); 
\end{align*}
by choosing $\delta < e^{\inf \{\lambda_n\}}$ (possible thanks to Step 2) and exploiting Step 4 (which implies that $|v_n|_{2^*_s}$ is bounded), we obtain that $|v_n|_2$ is bounded as well. 

\smallskip
\textbf{Step 6.} \emph{$(v_n)_n$ converges strongly in $H^s(\R^N)$ to a nontrivial function}.

In this step, we rely on the superlinearity of the perturbed problem.
From Steps 2--5, there exist $\lambda \in \R$ and $v \in H^s_{\textup{rad}}(\R^N)$ such that, up to a subsequence, $\lambda_n \to \lambda$, $v_n \rightharpoonup v$ in $H^s_{\textup{rad}}(\R^N)$, $v_n \to v$ in $L^q(\R^N)$ for $q \in (2,2^*_s)$, and $v_n(x) \to v(x)$ for a.e. $x \in \R^N$. In particular, from \hyperref[g0]{\rm{(g0)}}--\hyperref[g3]{\rm{(g3)}},
\[
\lim_n \intRN g_\eps(v_n) v \, \dx = \intRN g_\eps(v) v \, \dx \quad \text{and} \quad \lim_n \intRN g_+(v_n) v_n \, \dx = \intRN g_+(v) v \, \dx;
\]
moreover, from \eqref{eq_dim_extra_prop2} we have
\[
\partial_u I^m_{\eps}(\lambda_n,v_n)[v] = 0 = \partial_u I^m_{\eps}(\lambda_n,v_n)[v_n].
\]
Thus, exploiting the weak convergence, Fatou's lemma (recall $g_-^\eps(t) t \ge 0$ for all $t \in \R$) and that $\lambda_n \to \lambda$, we obtain
\begin{align*}
\MoveEqLeft|(-\Delta)^{s/2} v|_2^2 + e^\lambda |v|_2^2 - \intRN g_{\eps}(v) v \, \dx\\
&= \lim_n 
\intRN (-\Delta)^{s/2} v_n (-\Delta)^{s/2} v \, \dx + e^{\lambda_n} \intRN v_n v \, \dx - 
\intRN g_{\eps}(v_n) v \, \dx
 = 0 \\
&= \lim_n 
 \intRN |(-\Delta)^{s/2} v_n|^2 \, \dx + e^{\lambda_n} \intRN v_n^2 \, \dx - 
 \intRN g_{\eps}(v_n) v_n \, \dx 
 \\
&= \lim_n \left( \intRN |(-\Delta)^{s/2} v_n|^2 \, \dx + e^{\lambda_n} \intRN v_n^2 \, \dx + \intRN g_-^\eps(v_n) v_n \, \dx \right) - \intRN g_+(v) v \, \dx \\
&\geq \liminf_n \left( \intRN |(-\Delta)^{s/2} v_n|^2 \, \dx + e^{\lambda_n} \intRN v_n^2 \, \dx \right) - \intRN g_{\eps}(v) v \, \dx \\
&\geq |(-\Delta)^{s/2} v|_2^2 + e^\lambda |v|_2^2 - \intRN g_{\eps}(v) v \, \dx,
\end{align*}
thus 
$$ \liminf_n 
\intRN |(-\Delta)^{s/2} v_n|^2 \, \dx + e^{\lambda_n} \intRN v_n^2 \, \dx 
 = |(-\Delta)^{s/2} v|_2^2 + e^\lambda |v|_2^2;$$ 
this means that, up to a subsequence, $v_n \to v$ strongly in $H^s(\R^N)$. Finally, for $\delta\in (0,\inf_n e^{\lambda_n})$ and a suitable $C_{\delta}>0$,
$$
\frac{1}{2^*_s}|(-\Delta)^{s/2}v_n|_2^2 = - \frac{e^{\lambda_n}}{2} |v_n|_2^2 +\int_{\R^N} G_{\eps}(v_n) \, \dx \leq \frac{\delta-e^{\lambda_n}}{2} |v_n|_2^2 + C_{\delta} |v_n|_{2^*_s}^{2^*_s} \leq C_{\delta}' |(-\Delta)^{s/2}v_n|_2^{2^*_s}.
$$
This implies $|(-\Delta)^{s/2}v_n|_2^2$ is bounded away from zero and hence, by the strong convergence, $v \nequiv 0$.

\smallskip
\textbf{Step 7.} \emph{Conclusions}.

By the strong convergence, we have $I^m_{\eps}(\lambda, 
v)=d_{\eps,+}^m$, i.e.,
$(\lambda, v)$ is a Poho\v{z}aev minimum in the product space. The fact that $(\mu_\varepsilon,u_\varepsilon) := (e^\lambda,v)$ is a solution to \eqref{e-pert} follows from Proposition \ref{pr-Poho_sol}. 
Finally, assume $g|_{(-\infty,0)} = 0$ or $g$ odd, which imply the same properties on $g_{\eps}$.
Then, in Step 1, 
we may assume that each $v_n$ is Schwarz-symmetric in 
light of Theorem \ref{t-bks} (ii), hence so is $v$. 
\QED

\subsection{Symmetry results} 
\label{subsec_symmetry_post}

\indent

Exploiting the perturbed $g_{\eps}$, we prove here Proposition \ref{p-sym_unpert} (i)--(iii).
We recall that $G$ does not map $H^s(\R^N)$ to $L^1(\R^N)$; with the intention of gaining a symmetry result arguing as in \cite[Proof of Theorem 4.1]{LoMa08}, we show the following lemma.

\begin{Lemma}\label{l-LopMar}
Assume \hyperref[g0]{\rm{(g0)}}--\hyperref[g2]{\rm{(g2)}}. For $v \in H^s(\R^N)$ and $x = (x_1,x') \in \R \times \R^{N-1}$, define
\[
v_1(x) :=
\begin{cases}
v(x_1,x') \quad \text{if } x_1 < 0\\
v(-x_1,x') \quad \text{if } x_1 \ge 0
\end{cases}
\quad \text{and} \quad
v_2(x) :=
\begin{cases}
v(-x_1,x') \quad \text{if } x_1 < 0\\
v(x_1,x') \quad \text{if } x_1 \ge 0.
\end{cases}
\]
Then $\intRN G(v_1) \, \dx + \intRN G(v_2) \, \dx = 2 \intRN G(v) \, \dx$.
\end{Lemma}
\claim Proof. It is clear that $\intRN G_+(v_1) \, \dx + \intRN G_+(v_2) \, \dx = 2 \intRN G_+(v) \, \dx$. Additionally, from the monotone convergence theorem,
\begin{align*}
\intRN G_-(v_1) \, \dx + \intRN G_-(v_2) \, \dx & = \lim_{\eps \to 0^+} \intRN G_-^\eps(v_1) \, \dx + \intRN G_-^\eps(v_2) \, \dx\\
& = 2 \lim_{\eps \to 0^+} \intRN G_-^\eps(v) \, \dx = 2 \intRN G_-(v) \, \dx,
\end{align*}
whence the statement.
\QED

\bigskip

\claim Proof of Proposition \ref{p-sym_unpert} (i)--(iii). 

(i) Let $(\mu, u) \in \mc{P}_{(\rho)}$ be a minimizer as in the statement. In view of Remark \ref{rem_2min} and Lemma \ref{l-BJM}, $u$ solves the minimization problem \eqref{eq_minBS}. Then, exploiting Lemma \ref{l-LopMar}, we can argue in a similar way to \cite[Proof of Theorem 4.1]{LoMa08}\footnote{We observe that we formally deal with \textit{Case A} therein but minimize over $H^s(\R^N)$.}. Likewise if $u$ minimizes $K$ over $\mc{S}_m$.

%
%

(ii) 
Let $(\mu, u) \in \mc{P}_{(\rho)}$ be a minimizer as in the statement, and let us denote by $u^*$ the Schwarz rearrangement of $u$; in particular, $\int_{\R^N} G_-(u) \, \dx < +\infty$. 
As a consequence, from \cite[Proposition B.1]{BKS17},
\begin{equation}\label{e-Schwarz}
\intRN G_-(u) \, \dx = \intRN G_-(u^*) \, \dx.
\end{equation}
Now, assume by contradiction that $u$ does not equal (up to a translation) its Schwarz rearrangement $u^* \ne 0$. Then, from \cite[Propositions B.1 and B.3]{BKS17},
\[
|(-\Delta)^{s/2} u^*|_2^2 < |(-\Delta)^{s/2} u|_2^2 \quad \text{and} \quad |u^*|_2 = |u|_2,
\]
which imply
\[
J_{\mu}(u^*) < J_{\mu}(u) \quad \text{and} \quad P_{\mu}(u^*) < P_{\mu}(u) = 0.
\]
From Remark \ref{r-all_mu}, we find $t^* > 0$ such that $P_{\mu}\bigl(u^*(\cdot/t^*)\bigr) = 0$ and
\begin{align*}
\inf_{\mc{P}_{\mu}} J_{\mu} & \le J_{\mu}\bigl(u^*(\cdot/t^*)\bigr) = \frac{s}{N} |(-\Delta)^{s/2} u^*(\cdot/t^*)|_2^2 < \frac{s}{N} |(-\Delta)^{s/2} u^*|_2^2\\
& < \frac{s}{N} |(-\Delta)^{s/2} u|_2^2 = J_{\mu}(u) = \inf_{\mc{P}_{\mu}} J_{\mu},
\end{align*}
a contradiction. Similarly for $L^2$-minima. 

(iii) We claim that $\intRN |g_-(u) u^-| \, \dx < +\infty$, where, we recall, $u^\pm = \max\{\pm u,0\}$. 
As a matter of fact, since $g_-(u) u^- = g_-(-u^-) u^- \le 0$, there holds
\[
|g_-(u) u^-| = -g_-(-u^-) u^- \le g_-(u^+) u^+ - g_-(-u^-) u^- = g_-(u) u^+ - g_-(u) u^- = g_-(u)u \in L^1(\R^N),
\]
which holds by the assumptions.
As a consequence, 
with the notation of Lemma \ref{lem_converg_Hs}, 
	for every $R>0$ there holds
	\begin{equation*}
		\intRN (-\Delta)^{s/2}u \ (-\Delta)^{s/2}(\varphi_R u^-) \, \dx + \mu \intRN u \varphi_R u^- \, \dx = \intRN g(u) \varphi_R u^- \, \dx,
	\end{equation*}
	and passing to the limit we obtain
\begin{equation}\label{e-ROFL}
\intRN (-\Delta)^{s/2} u \ (-\Delta)^{s/2} u^- \, \dx + \mu \intRN u u^- \, \dx = \intRN g(u) u^- \, \dx.
\end{equation}
Finally, using \eqref{e-ROFL} and Lemma \ref{l-assume+},
\begin{align*}
	0 & \le 
	\intRN \left(g(-u^-)- \mu (-u^-) \right) u^- \, \dx = \intRN g(u) u^- - \mu u u^- \, \dx
	 = \intRN (-\Delta)^{s/2}u (-\Delta)^{s/2}u^- 
	\, \dx\\
	& \le -|(-\Delta)^{s/2}u^-|_2^2 
	 \le 0,
\end{align*}
and we conclude. \QED
%

\section{The 
sublinear problem}\label{s:limit}

\subsection{$\eps$-independent estimates}


\indent

Let $\mu_{\eps} = e^{\lambda_{\eps}}$, $u_{\eps}$, $m_{\eps}$, and $d_{\eps,+}^m$ 
be the quantities found in Theorem \ref{t-pert} and $m_0$ be as in \eqref{eq_meps_m0}. 
From now on, we are interested in the dependence on $\eps$; in particular, we will ensure that some objects related to problem \eqref{e-pert} enjoy some uniformity in $\eps > 0$.
We start by showing that $
(d_{\eps,+}^m)_{\eps \in (0,1]}$ is bounded and bounded away from zero.

\begin{Lemma}\label{lem_uniform_energ_bound}
Let $m>m_0$. There exist $\zeta_m, \beta_m \in \R$ 
such that
$$-\infty < \zeta_m \leq d_{\eps,+}^m 
\leq \beta_m<0 \quad \hbox{ for every $\eps \in (0,1]$}.$$
\end{Lemma}

\claim Proof.
By Lemma \ref{l-bddB} and the monotonicity of $\eps \mapsto a_{\eps}(\lambda)$ and $\eps \mapsto \overline{\lambda}_{\eps}$ -- cf. \eqref{eq_comp_aeps_a0} and \eqref{e-mu0}, we have $d_{\eps,+}^m=b_{\eps}^m \geq b_1^m=: \zeta_m$. 
Again by Lemma \ref{l-bddB} 
we have, for every $\lambda < 
\overline{\lambda}_0 \leq \overline{\lambda}_{\eps}$,
$$d_{\eps,+}^m 
\leq \frac{e^{\lambda}}{2} \left( \frac{a_{\eps}(\lambda)}{e^{\lambda}/2} -m\right) \leq \frac{e^{\lambda}}{2} \left( \frac{a(\lambda)}{e^{\lambda}/2} -m\right).$$
Let now $\delta := \frac{m-m_0}{2}>0$. By the definition of infimum $m_0 = \inf_{\lambda < \overline{\lambda}_0} \frac{a(\lambda)}{e^\lambda/2}$ there exists $\lambda_{\delta} < \overline{\lambda}_0$ such that
$$d_{\eps,+}^m 
\leq \frac{e^{\lambda_{\delta}}}{2} \left( \frac{a(\lambda_{\delta})}{e^{\lambda_{\delta}}/2} -m\right) \leq \frac{e^{\lambda_{\delta}}}{2} \left( m_0+\delta-m\right) = -\frac{e^{\lambda_{\delta}}}{4} (m-m_0) =: \beta_m <0. 
\QED
$$

\bigskip



We show now some estimates on the solutions $(\lambda_{\eps}, u_{\eps})$. Here we use the $L^2$-subcritical growth \hyperref[g3]{\rm{(g3)}}.


\begin{Lemma}\label{lem_bounded_sequence_eps}
Let $m>m_0$. Then $(\lambda_{\eps}, u_{\eps})_{\eps \in (0,1]}$ is bounded in $\R \times H^s_{\textup{rad}}(\R^N)$. 
\end{Lemma}

\claim Proof.
Observe first that, since $(\lambda_{\eps}, u_{\eps})$ satisfies the Poho\v{z}aev identity and by Lemma \ref{lem_uniform_energ_bound}, we have
\begin{equation*}
\frac{m}{2} e^{\lambda_{\eps}} = \frac{1}{2^*_s} P_\eps(\lambda_{\eps},u_{\eps}) - I^m_{\eps}(\lambda_{\eps},u_{\eps}) + \frac{s}{N} |(-\Delta)^{s/2} u_{\eps}|_2^2 \ge - d_{\eps,+}^m 
\geq - \beta_m
\end{equation*}
which gives a bound below on $(\lambda_{\eps})_{\eps \in (0,1]}$. 
Next, from \hyperref[g0]{\rm{(g0)}}, \hyperref[g1]{\rm{(g1)}}, and \hyperref[g3]{\rm{(g3)}}, in a similar way to \eqref{e-blue}, for every $\delta>0$ there exists $c_\delta>0$ such that for all $t \in \R$
\[
g_+(t)t \le 
c_{\delta} t^2 + \delta |t|^{\bar{p}+1}.
\]
Then, from \eqref{eq_GN_ineq} 
and the fact that $(\lambda_{\eps}, u_{\eps})$ is a critical point of $I_\eps^m$, we have
\begin{align*}
0 &= |(-\Delta)^{s/2} u_{\eps}|_2^2 + e^{\lambda_{\eps}} |u_{\eps}|_2^2 + \intRN g_-^\eps(u_{\eps}) u_{\eps} \, \dx - \intRN g_+(u_{\eps}) u_{\eps} \, \dx\\
&\ge |(-\Delta)^{s/2} u_{\eps}|_2^2 + \left(e^{\lambda_{\eps}} - 
c_\delta\right) |u_{\eps}|_2^2 + \intRN g_-^\eps(u_{\eps}) u_{\eps} \, \dx - \delta |u_{\eps}|_{\bar{p}+1}^{\bar{p}+1} 
\\
&\ge |(-\Delta)^{s/2} u_{\eps}|_2^2 + \left(e^{\lambda_{\eps}} -c_\delta\right) |u_{\eps}|_2^2 - \delta C_\textup{GN} |(-\Delta)^{s/2} u_{\eps}|_2^2 |u_{\eps}|_2^{4s/N} \\
&= \left(1- \delta C_\textup{GN} m^{2s/N} \right) |(-\Delta)^{s/2} u_{\eps}|_2^2 + \left(e^{\lambda_{\eps}} - c_\delta\right) m,
\end{align*}
and we conclude by taking $\delta < (C_\textup{GN} m^{2s/N})^{-1}$.
\QED

%

\subsection{Passage to the limit}
\label{subsec_passage_lim}

\indent

We want to pass to the limit with the solutions $(\lambda_{\eps}, u_{\eps})$ found in Theorem \ref{t-pert}. Observe that, from \eqref{eq_compar_meps_m0}, if $m > m_0$, then $m > m_\eps$ for every $\eps \in (0,1]$. Recall that, by
Lemma \ref{lem_uniform_energ_bound} and Proposition \ref{prop_uguag_not_sym}, we have
$$-\infty < 
\zeta_m \leq d_{\eps,+}^m = I_\varepsilon^m(\lambda_\varepsilon,u_\varepsilon) 
= \kappa^m_{\eps} = \ell^m_{\eps} = K_\varepsilon(u_\varepsilon) \leq \beta_m <0;$$
moreover, by Lemma \ref{lem_bounded_sequence_eps}, $(\lambda_{\eps}, u_{\eps})_{\eps \in (0,1]}$ is bounded in $\R \times H^s_{\textup{rad}}(\R^N)$. Therefore, up to a subsequence,
%
%
%
%
%
%
%
\begin{equation}\label{eq_converg_u_eps}
\begin{cases}
u_{\eps} \wto u_0 & \quad \hbox{in $H^s_{\textup{rad}}(\R^N)$}, \\ 
\lambda_{\eps} \to \lambda_0 & \quad \hbox{in $\R$}.
\end{cases}
\end{equation}
In addition, we have
$$|u_0|_2^2 \leq m$$
by the weak convergence. 
%
We write $\mu_{\eps}:= e^{\lambda_{\eps}}$, $\mu_0:= e^{\lambda_0}$.



\smallskip

We want to show that $(\mu_0, u_0) \in (0,+\infty)\times H^s_{\textup{rad}}(\R^N)$ is the desired solution of Theorem \ref{t-main}. To this aim, the first step is to prove the pointwise convergence of $G_\eps(u_\eps)$ to $G(u_0)$.

\begin{Lemma}\label{lem_conv_ae} 
Let $u_0 \in H^s_\textup{rad}(\R^N)$ as in \eqref{eq_converg_u_eps}. Then, $g_\eps(u_\eps(x)) \to g(u_0(x))$ and $G_\eps(u_\eps(x)) \to G(u_0(x))$ as $\varepsilon \to 0^+$ for a.e. $x \in \R^N$. 
\end{Lemma}

%
%

\claim Proof.
Let $x \in \R^N$ such that $u_{\eps}(x) \to u_0(x) \in \R$.
We prove first the claim on $g_{\eps}$.
We write
$$(0,1] = \{ \eps \in (0,1] \mid |u_{\eps}(x)|>\eps\} \cup \{ \eps \in (0,1] \mid |u_{\eps}(x)|\leq \eps\} =: I^x_1 \cup I^x_2.$$
Clearly, $0 \in \overline{I^x_1} \cup \overline{I^x_2}$. 
%
 If $0 \in \overline{I^x_1}$ then, observed that $g_{\eps}(u_{\eps}(x))=g(u_{\eps}(x))$ for every $\eps \in I^x_1$, we get
$$\lim_{\eps \to 0^+, \; \eps \in I^x_1} g_{\eps}(u_{\eps}(x)) = \lim_{\eps \to 0^+, \; \eps \in I^x_1} g(u_{\eps}(x)) = g(u_0(x)).$$
%
If $0 \in \overline{I^x_2}$, let $\delta>0$. By the assumptions, we know that there exists $\omega_{\delta}>0$ such that $|g(t)|\leq \delta$ for every $|t|\leq \omega_{\delta}$. Observed that 
(recall \eqref{eq_intr_g+}-\eqref{eq_intr_g-}) 
\begin{align*}
|g_{\eps}(u_{\eps}(x))| &\leq |g_+(u_{\eps}(x))| + |g_-^{\eps}(u_{\eps}(x))| \leq |g_+(u_{\eps}(x))| + |g_-(u_{\eps}(x))| \\
&= g^+(u_{\eps}(x)) + g^-(u_{\eps}(x)) 
= |g(u_{\eps}(x))|\leq \delta
\end{align*}
for every $\eps \in I^x_2 \cap (0, \omega_{\delta}]$, we obtain
$$\limsup_{\eps \to 0^+, \; \eps \in I^x_2} |g_{\eps}(u_{\eps}(x)) | \leq \delta;$$
by the arbitrariness of $\delta$, we achieve
$$\lim_{\eps \to 0^+, \; \eps \in I^x_2} g_{\eps}(u_{\eps}(x)) =0;$$
on the other hand, by definition of $I^x_2$, we have
$$\lim_{\eps \to 0^+, \; \eps \in I^x_2} |u_{\eps}(x)| \leq \lim_{\eps \to 0^+, \; \eps \in I^x_2} \eps = 0$$
and thus $u_0(x) = \lim_{\eps \to 0^+, \; \eps \in I^x_2} u_{\eps}(x)=0$; 
since $g(u_0(x))=g(0)=0$, we have reached
$$\lim_{\eps \to 0^+, \; \eps \in I^x_2} g_{\eps}(u_{\eps}(x)) =g(u_0(x)),$$
which is the claim.
%
We prove now the claim on $G_{\eps}$.
Obviously, $G_+(u_\eps(x)) \to G_+(u_0(x))$; 
noticed that, in a similar way as before, we can show $g_-^{\eps}(u_{\eps}
(x)t) \to g_- 
(u_0(x)t)$ for each $t \in [0,1]$, and observed that
$$
|g_-^\eps(u_\eps(x) t)| \le |g_-(u_\eps(x)t)| \lesssim 
 1 + |u_{\eps}(x)|^{2^*_s} \leq 2 + |u_0(x)|^{2^*_s} 
$$
for $0 < \eps \ll 1$ and any $t \in [0,1]$, 
 we can use the dominated convergence theorem and obtain
$$
G_-^{\eps}(u_{\eps}(x)) = u_{\eps}(x) \int_0^1 g_-^{\eps}(u_{\eps}(x)t) \, \mathrm{d}t \to u_0(x) \int_0^1 g_-(u_0(x)t) \, \mathrm{d}t = G_-(u_0(x)). \QED
$$

\medskip

We observe that the limit function obtained in \eqref{eq_converg_u_eps}, at the moment, could be trivial. In what follows, we show that this is not the case and, moreover, $u_0$ is a solution enjoying minimality properties.

\begin{Proposition}\label{prop_lm_km_em}
Let $(\mu_0,u_0) \in (0,+\infty) \times H^s_\textup{rad}(\R^N)$ as in \eqref{eq_converg_u_eps}. Then $u_0 \in \mc{S}_m$, $(\mu_0,u_0) \in \mc{P}$, it is a solution to \eqref{e-main}, and
\begin{equation*}
K(u_0) = I^m(\mu_0,u_0) = \ell^m = \kappa^m = d^m_+ = d^m < 0.
\end{equation*}
Moreover, $u_{\eps} \to u_0$ in $H^s(\R^N)$, $\intRN G_-^{\eps} (u_{\eps}) \, \dx \to \intRN G_-(u_0) \, \dx<\infty$, and $d^m_{\eps,+}\to d^m_+$ as $\varepsilon \to 0^+$.
\end{Proposition}

\claim Proof.
We show first that $u_0$ minimizes $K$ over the $L^2$-ball.
Since $G_{\eps} \geq G$, we have
$$ \ell_{\eps}^m 
= \inf_{\mc{D}_m^\textup{rad}} K_{\eps} \leq \inf_{\mc{D}_m^\textup{rad}} K \leq K(u_0).$$
On the other hand, by exploiting that $(-\Delta)^{s/2} u_{\eps} \wto (-\Delta)^{s/2} u_0$ 
in $L^2(\R^N)$, the strong convergence in $L^1(\R^N)$ of 
$G_+(u_{\eps})$, and Fatou's lemma for $G_-^{\eps}(u_{\eps})$ (here we use Lemma \ref{lem_conv_ae}),
\begin{align*}
K(u_0) &= \frac{1}{2} |(-\Delta)^{s/2} u_0|_2^2 
- \intRN G_+(u_0) \, \dx + \intRN G_-(u_0) \, \dx \\
& \leq \liminf_{\eps \to 0^+} \left( \frac{1}{2} |(-\Delta)^{s/2} u_{\eps} |_2^2\right) 
- \lim_{\eps \to 0^+} \left( \intRN G_+(u_{\eps}) \, \dx \right) + \liminf_{\eps \to 0^+} \left( \intRN G_-^{\eps}(u_{\eps}) \, \dx \right) \\
& \leq \liminf_{\eps \to 0^+} K_{\eps}(u_{\eps}) = \liminf_{\eps \to 0^+} \ell_{\eps}^m \le \limsup_{\eps \to 0^+} \ell_{\eps}^m \le \inf_{\mc{D}_m^\textup{rad}} K \le K(u_0). 
\end{align*}
This yields that $u_0$ is a minimum over the $L^2$-ball and the following holds as $\eps \to 0^+$: $\ell_{\eps}^m \to \ell^m$, $\intRN G_-^{\eps} (u_{\eps}) \, \dx \to \intRN G_-(u_0) \, \dx<\infty$, and $(-\Delta)^{s/2} u_{\eps} \to (-\Delta)^{s/2} u_0$ in $L^2(\R^N)$.
By the fact that $ \ell_{\eps}^m \leq \beta_m < 0$, 
we obtain $K(u_0) \leq \beta_m 
< 0$; in particular, $u_0 \neq 
0$. Moreover, by 
Proposition \ref{pr:PohoMin} (iv) we have that $u_0 \in \mc{S}_m$ which implies $u_{\eps} \to u_0$ in $H^s(\R^N)$, and hence $P(\mu_0,u_0)=0$. Moreover, by Proposition \ref{pr:PohoMin} (i) there exists $\tilde{\mu}$ such that $(\tilde{\mu},u_0)$ is a solution of \eqref{eq_general_mu_fixed} and it is a Poho\v{z}aev minimum over $ \mc{P}_{\tilde{\mu}}$; since $P(\mu_0,u_0)=0=P(\tilde{\mu},u_0)$ and $u_0 \neq 0$, it must be $\tilde{\mu}=\mu_0$.
Finally by Proposition \ref{prop_uguag_not_sym} we have $\kappa^m = \ell^m = d^m_+ = d^m$.
\QED

\bigskip

To show the Nehari identity, we need the following lemma.

\begin{Lemma}\label{lem_stima_nehar}
The solution $u_0$ 
given in \eqref{eq_converg_u_eps} satisfies $\intRN g_-(u_0)u_0 \, \dx < \infty$.
\end{Lemma}

\claim Proof.
Let us recall that, as $\eps \to 0^+$, 
$\mu_{\eps}\to \mu_0$, 
$u_{\eps} \to u_0$ in $H^s(\R^N)$ and a.e. in $\R^N$, and 
$g_\eps(u_\eps) \to g(u_0)$ a.e. in $\R^N$. Additionally, $(\mu_\eps,u_\eps)$ satisfies the Nehari identity
\[
|(-\Delta)^{s/2} u_\eps|_2^2 + \mu_\eps |u_{\eps}
|_2^2 = \intRN g_+(u_\eps) u_{\eps}- g_-^\eps(u_\eps)u_{\eps} \, \dx.
\]
Since $g_-(t)t \ge 0$ and $\varphi_\eps(t) \ge 0$ for all $t \in \R$, we also have that $g_-^\eps(t) t \ge 0$ for all $t \in \R$, hence $u_0$ satisfies the following ``Nehari inequality''
\begin{align*}
|(-\Delta)^{s/2} u_0|_2^2 + \mu_0 |u_0|_2^2 + \intRN g_-(u_0) u_0 \, \dx & \le \lim_{\eps \to 0^+} |(-\Delta)^{s/2} u_\eps|_2^2 + \mu_\eps |u_\eps|_2^2 + \intRN g_-^\eps(u_\eps) u_\eps \, \dx\\
& = \lim_{\eps \to 0^+} \intRN g_+(u_\eps) u_\eps \, \dx = \intRN g_+(u_0) u_0 \, \dx
\end{align*}
thanks to \hyperref[g1]{\rm{(g1)}}, \hyperref[g3]{\rm{(g3)}}, and Fatou's lemma. In particular, this gives the claim.
\QED

\bigskip

\medskip

\claim Proof of Theorem \ref{t-main}.
The theorem is a consequence of 
Lemma \ref{lem_stima_nehar}, Proposition \ref{p-Nehari}, Proposition \ref{prop_lm_km_em}, Proposition \ref{prop_intr_legendre} and Proposition \ref{prop_JJLU}.
\QED
%

\bigskip


\claim Proof of Proposition \ref{p-sym_unpert} (iv).
From Theorem \ref{t-pert} and \eqref{eq_compar_meps_m0}, for every $\varepsilon \in (0,1]$ there exists $(\lambda_\varepsilon,u_\varepsilon) \in \mc{P}_{\varepsilon,+}$ Schwarz-symmetric such that $I^m_\varepsilon(\lambda_\varepsilon,u_\varepsilon) = \inf_{\mc{P}_{\varepsilon,+}} I^m_\varepsilon = \inf_{\mc{P}_{\varepsilon,+}^\textup{rad}} I^m_\varepsilon$. Then, repeating all the steps in Section \ref{s:limit}, not only do we re-obtain Theorem \ref{t-main}, but we also get the same outcome where $\kappa^m$, $\ell^m$, $a(\mu)$, and $d^m$ are replaced, respectively, with $\inf_{\mc{S}_m} K$, $\inf_{\mc{D}_m} K$, $\inf_{\mc{P}_\mu} J_\mu$, and $\inf_{\mc{P}} I^m$. 
\QED

\section{Further existence results}\label{s-FER}

\subsection{Existence for all masses}

\indent

In Theorem \ref{t-main}, we obtained the existence of a solution for large masses. When we rule out the possibility of $g$ to be superlinear at the origin (or positive nearby), which clearly mimics the logarithmic case, we can achieve existence also for all masses. Obviously, we already have existence for all masses if $m_0 = 0$, which is why we assume $m_0 > 0$ in Theorem \ref{thm_small_mass} below.
%
%
%
%

%

\begin{Lemma}\label{l-mitilde}
Under assumptions \hyperref[g0]{\rm{(g0)}}--\hyperref[g4]{\rm{(g4)}}, for every $\omega \geq 0$ there exists $\widetilde{\mu}_{\omega} \in (0,+\infty)$ such that
\begin{equation}\label{e-bdd_interval}
\inf_{0 < \mu < \overline{\mu}_0} \frac{a(\mu)}{\mu + \omega'} =  \inf_{0 < \mu < \widetilde{\mu}_{\omega}} \frac{a(\mu)}{\mu + \omega'}
\end{equation}
for every $\omega' \in [0,\omega]$. In particular, $m_0 =  \inf_{0 < \mu < \widetilde{\mu}_0} \frac{a(\mu)}{\mu/2}.$
Moreover, 
$$ \frac{\widetilde{\mu}_{\omega}}{\omega} \to 0 \quad \hbox{as $\omega \to +\infty$}.$$
\end{Lemma}
\claim Proof.
If $\overline{\mu}_0 < +\infty$, we take $\widetilde{\mu}_{\omega} := \overline{\mu}_0$. 
Otherwise, we set $M:= \inf_{\mu \in (0,+\infty)} \frac{a(\mu)}{\mu}$ and
\begin{equation*}
\widetilde{\mu}_{\omega} := \inf\left\{ \widetilde{\mu}>0  \;\middle|\;
 \frac{a(\mu)}{\mu+\omega} \geq M+1 \; \hbox{for each $\mu \geq \widetilde{\mu}$} \right\}.
\end{equation*}
We observe that $\widetilde{\mu}_{\omega} < +\infty$ from Lemma \ref{l-bddB}, while $\widetilde{\mu}_{\omega} > 0$ from its definition and that of $M$. Additionally, from the continuity of $a \colon (0,+\infty) \to [0,+\infty)$ (see Proposition \ref{p-a_monot}), we have
\begin{equation}\label{e-aM+1}
\frac{a(\widetilde{\mu}_{\omega})}{\widetilde{\mu}_{\omega} + \omega} = M+1.
\end{equation}
Consequently, if $\nu \ge \widetilde{\mu}_{\omega}$ and $\omega' \in [0,\omega]$, then
\[
\frac{a(\nu)}{\nu + \omega'} \ge M+1 > \inf_{0 < \mu < +\infty} \frac{a(\mu)}{\mu + \omega'},
\]
whence \eqref{e-bdd_interval}.
We move to the last claim; we can assume $\widetilde{\mu}_{\omega}\to +\infty$ as $\omega \to +\infty$ (otherwise, the claim is obvious). 
Since \eqref{e-aM+1} can be rewritten as
$$\frac{\omega}{\widetilde{\mu}_{\omega}} = \frac{1}{M+1} \frac{a(\widetilde{\mu}_{\omega})}{\widetilde{\mu}_{\omega}} - 1,$$
the claim follows from 
Lemma \ref{l-bddB}.
\QED


\begin{Theorem}\label{thm_small_mass}
Assume \hyperref[g0]{\rm{(g0)}}--\hyperref[g4]{\rm{(g4)}}, $m_0 > 0$, and
$$\limsup_{t \to 0} \frac{g(t)}{t} =: - \overline{\omega} \in [-\infty, 0).$$
Then, there exists $m_1=m_1(\overline{\omega}) \in [0, m_0)$ such that, for every $m > m_1$, there exists a 
solution $(\mu_0, u_0) \in \R \times \mc{S}_m^\textup{rad}$ to \eqref{e-main}, which satisfies the Nehari and Poho\v{z}aev identities, and such that
$$I^m(\mu_{0}, u_0) = K(u_0) = a(\mu_0) - \frac{\mu_0}{2} m= \kappa^m 
= d^m = \inf_{-\infty < \nu < \overline{\mu}_0} \left(a(\nu)-\frac{\nu}{2}m\right);$$
%
%
if $\kappa^m \leq 0$ 
then $\mu_0>0$. 
Finally, if $\overline{\omega}=\infty$, then $m_1=0$.
\end{Theorem}


\claim Proof.
Let $\omega \in (0,\overline{\omega})$ and define, for $t \in \R$,
\begin{equation}\label{e-h_omega}
h_{\omega}(t):= g(t) + \omega t, \quad
H_{\omega}(t):= G(t) + \frac{\omega}{2} t^2.
\end{equation}
We observe that
\begin{itemize}
\item $h_{\omega}$ is continuous and $h_{\omega}(0)=0$,
\item $\limsup_{|t|\to +\infty} \frac{|h_{\omega}(t)|}{|t|^{2^*_s-1}} < \infty$,
	\item there exists $t_0 \neq 0$ such that $H_{\omega}(t_0) >0$.
\end{itemize}
Moreover, it is straightforward to check that $ \limsup_{t \to 0} \frac{h_{\omega}(t)}{t} < 0$ and $\limsup_{|t| \to +\infty} \frac{h_{\omega}(t)}{|t|^{\bar{p}-1}t} \leq 0$, thus
$h_{\omega}$ satisfies the assumptions of Theorem \ref{t-main}. In particular, for
\begin{equation*}
m>m_{\omega}:= \inf_{0<\mu < \overline{\mu}_{\omega}} \frac{a_{\omega}(\mu)}{\mu/2}
\end{equation*}
there exists an $L^2$-minimum $u_{\omega} \in \mc{S}_m^\textup{rad}$ with Lagrange multiplier $\mu_{\omega}>0$, i.e.,
\begin{equation}\label{eq_minim_omega_smallm}
K_{\omega}(u_{\omega}) = \inf_{\mc{D}_m^\textup{rad}} K_{\omega} = \inf_{\mc{S}_m^\textup{rad}} K_{\omega} =: \kappa^m_{\omega}<0;
\end{equation}
here, with the convention that $+\infty + t = +\infty$ for every $t \in \R$,
$$\overline{\mu}_{\omega}:= \sup_{t \neq 0} \frac{H_{\omega}(t)}{t^2/2} = \overline{\mu}_0 + \omega, \quad K_{\omega}(u):= \frac{1}{2} |(-\Delta)^{s/2}u|_2^2 - \intRN H_{\omega}(u) \, \dx,$$
$$a_{\omega}(\mu):=\inf_{\mc{P}_{\mu,\omega}^{\textup{rad}}} J_{\mu,\omega}, \quad J_{\mu,\omega}(u):= \frac{1}{2} |(-\Delta)^{s/2}u|_2^2 + \frac{\mu}{2} |u|_2^2- \intRN H_{\omega}(u) \, \dx,$$
$$P_{\mu,\omega}(u) 
:= \intRN |(-\Delta)^{s/2} u|^2 \, \dx + 2^*_s \intRN \frac{\mu}{2} u^2 - H_{\omega}(u) \, \dx,$$
$$\mc{P}_{\mu,\omega} := \left\{u \in H^s(\R^N) \setminus \{0\} \mid P_{\mu,\omega}(u) = 0\right\}, \quad \mc{P}_{\mu,\omega}^{\textup{rad}} := \mc{P}_{\mu,\omega} \cap H^s_{\textup{rad}}(\R^N).$$
We notice that $(\mu_{\omega},u_{\omega})$ satisfies the Nehari and Poho\v{z}aev identities with respect to $H_{\omega}$.
Since 
$K_{\omega}(u) = K(u) - \frac{\omega}{2} |u|_2^2$,
we obtain 
$$ K(u_{\omega}) - \frac{\omega}{2}m = K(u_{\omega}) - \frac{\omega}{2}|u_{\omega}|_2^2 = K_{\omega}(u_{\omega}) = \inf_{\mc{S}_m^\textup{rad}} K_{\omega} = \kappa^m - \frac{\omega}{2} m,$$
whence $K(u_{\omega}) = \kappa^m$, 
which means that $u_{\omega}$ is an $L^2$-minimum for the original problem. 
Moreover,
\begin{equation}\label{eq_energ_nosign}
\kappa^m=\kappa^m_{\omega} + \frac{\omega}{2} m.
\end{equation}
In addition, $u_{\omega}$ is a 
solution to the Euler--Lagrange equation corresponding to $K_{\omega}$, i.e., 
$$(-\Delta)^{s} 
u_{\omega} + \mu_{\omega} u_{\omega} = h_{\omega}(u_{\omega}) \quad \hbox{ in $\R^N$},$$
which means
that $u_0 := u_{\omega}$ is a 
solution 
to \eqref{e-main} with Lagrange multiplier
\begin{equation}\label{eq_molt_lagr_nosign}
\mu_0:= \mu_{\omega} - \omega
\end{equation}
and 
$(\mu_0,u_0)$ satisfies the Nehari and Poho\v{z}aev identities with respect to $G$.
%
%
Additionally, for any $\nu>\omega$, we have
$J_{\nu,\omega}(u) = J_{\nu-\omega}(u) $ and $P_{\nu,\omega}(u) = P_{\nu-\omega}(u) $, 
therefore
\begin{equation}\label{eq_momega}
\frac{1}{2}m_{\omega}
\leq \inf_{\omega<\nu <\overline{\mu}_0+\omega} \frac{a_{\omega}(\nu)}{\nu} 
= \inf_{0<\nu-\omega <\overline{\mu}_0} \frac{a(\nu-\omega)}{\nu}
= \inf_{0<\nu <\overline{\mu}_0} \frac{a(\nu)}{\nu+\omega} 
= \inf_{0<\nu <\overline{\mu}_0} \frac{a(\nu)}{\nu} \frac{\nu}{\nu + \omega}.
\end{equation}
Let $\widetilde{\mu} = \widetilde{\mu}_{\omega} \in (0,+\infty)$ be the value given in Lemma \ref{l-mitilde}.
Since the function $\nu \in (0,+\infty) \mapsto \frac{\nu}{\nu + \omega} \in (0,+\infty)$ is increasing, from \eqref{eq_momega} we obtain
\begin{equation*}
\frac12 m_\omega \le \inf_{0 < \nu < \widetilde{\mu}} \frac{a(\nu)}{\nu} \frac{\nu}{\nu + \omega} \le \frac{\widetilde{\mu}}{\widetilde{\mu} + \omega} \inf_{0 < \nu < \widetilde{\mu}} \frac{a(\nu)}{\nu} = \frac{\widetilde{\mu}}{\widetilde{\mu} + \omega} \inf_{0 < \nu < \overline{\mu}_0} \frac{a(\nu)}{\nu} = \frac{\widetilde{\mu}}{\widetilde{\mu} + \omega} \frac{m_0}{2} < \frac12 m_0.
\end{equation*}
In addition, if $\overline{\omega}=\infty$, taking $\omega$ arbitrarily large we obtain 
\begin{equation}\label{eq_momega_zero}
m_{\omega} \to 0 \quad \hbox{as $\omega \to +\infty$}
\end{equation}
again from Lemma \ref{l-mitilde}. As a consequence, for any $m>0$, there exists $m_{\omega}\in [0, m)$ and thus an $L^2$-minimum $u_{\omega}$. We conclude by setting
$$
m_1 := \inf_{0 < \omega < \overline{\omega}} m_{\omega}; 
$$
noticed that $\omega \mapsto m_\omega$ is nonincreasing (this is a direct consequence of $J_{\mu,\omega} \leq J_{\mu}$, $P_{\mu, \omega} \leq P_{\mu}$, $\mu_{\omega} \geq \mu_0$, and an argument similar to \eqref{eq_comp_aeps_a0}--\eqref{eq_compar_meps_m0}), we further observe that $m_1 = \lim_{\omega \to \overline{\omega}^-} m_{\omega}$. 
%
The sign of $\mu_0$ is a simple consequence of the Poho\v{z}aev identity.
Finally, the fact that
\[
K(u_0) = I^m(\mu_0, u_0) = d^m = \inf_{-\infty < \nu < \overline{\mu}_0} \left(a(\nu)-\frac{\nu}{2}m\right) = a(\mu_0) - \frac{\mu_0}{2} m
\]
follows from Propositions \ref{prop_uguag_not_sym}, \ref{prop_intr_legendre}, and \ref{prop_JJLU}. \QED

\begin{Remark}
By definition, $m_1$ is nonincreasing with respect to $\overline{\omega}$. Moreover, in the case of small masses, Theorem \ref{thm_small_mass} gives no information on the sign of the Lagrange multiplier and of the minimal energy -- see indeed \eqref{eq_energ_nosign} and \eqref{eq_molt_lagr_nosign}.
\end{Remark}


\claim Proof of Theorem \ref{t-small-mass}.
It is a consequence of Theorem \ref{thm_small_mass}. 
\QED

\subsection{Existence under relaxed assumptions}

\indent

We deal now with Proposition \ref{prop_esist_no_groundstate}.
We thus search for a normalized solution which is nonnegative, possibly considering also the case 
$$\limsup_{|t| \to +\infty} \frac{g_-(t)}{|t|^{2_s^*-1}} = +\infty.$$
The idea is to modify the nonlinearity $g$. We begin with the following result.
\begin{Lemma}\label{lem_apriori_bound}
Let \hyperref[g0]{\rm{(g0)}}--\hyperref[g2]{\rm{(g2)}} hold. Assume $g(t)=0$ for $|t| \geq t_1$, and let $u$ be a nonnegative solution of \eqref{eq_general_mu_fixed} for some $\mu \geq 0$. Then $u\in L^{\infty}(\R^N)$ and $|u|_{\infty} \leq t_1$.
\end{Lemma}


\claim Proof.
Let $v:=(u-t_1)^+$. It is easy to check, via explicit computations, that $v \in H^s(\R^N)$. 
Let $\varphi_R$ be as in Lemma \ref{lem_converg_Hs} and set $v_R := v 
\varphi_R$. 
Then $v_R \in H^s(\R^N) \cap L^1(\R^N)$ and hence
$$ 
\intRN (-\Delta)^{s/2} u (-\Delta)^{s/2} v_R \, \dx + \mu 
\intRN u v_R \, \dx = \intRN g(u) v_R \, \dx.$$
We notice that $g(u) v_R \equiv 0$ by the assumptions. 
Since $u \geq0 $ and $\mu\geq 0$, we have $\mu u v_R \geq 0$ and thus
$$ \intRN (-\Delta)^{s/2} u (-\Delta)^{s/2} v_R \, \dx \leq 0.$$
By Lemma \ref{lem_converg_Hs} we have
$$ \intRN (-\Delta)^{s/2} u (-\Delta)^{s/2} (u-t_1)^+ \, \dx \leq 0.$$
On the other hand, using the inequality $(a-b)(a^+-b^+) \geq |a^+-b^+|^2$ for every $a,b \in \R$, 
\begin{align*}
\MoveEqLeft \intRN (-\Delta)^{s/2} u (-\Delta)^{s/2} (u-t_1)^+ \, \dx \\
&=\frac{C_{N,s}}{2} \intRN \intRN \frac{\big((u(x)-t_1)-(u(y)-t_1)\big)\big((u(x)-t_1)^+ - (u(y)-t_1)^+\big)}{|x-y|^{N+2s}} \\
&\geq
\frac{C_{N,s}}{2} \intRN \intRN \frac{\big|(u(x)-t_1)^+ - (u(y)-t_1)^+\big|^2}{|x-y|^{N+2s}} = |(-\Delta)^{s/2} (u-t_1)^+|_2^2.
\end{align*}
Joining the two relations found we obtain $(u-t_1)^+=0$, i.e., $u \leq t_1$. 
\QED

\bigskip

\claim Proof of Proposition \ref{prop_esist_no_groundstate}.
Let us assume that $t_0$, defined in \hyperref[g4]{\rm{(g4)}}, is positive (if it is negative, then we argue likewise). We can assume $g(t_0)>0$ as well. Let $t_1$ be the first zero of $g$ larger than $t_0$ (if $g(t) > 0$ for all $t > t_0$, we set $t_1 = +\infty$). 
Then, we set
$$\overline{g}(t):= g(t) \chi_{(0, t_1)}(t).$$
Clearly, $\overline{g}$ satisfies \hyperref[g0]{\rm{(g0)}}--\hyperref[g4]{\rm{(g4)}}, 
thus Proposition \ref{p-sym_unpert} (iv) 
applies, 
and we have a Schwarz-symmetric solution $(\mu_0,u_0) \in (0,+\infty) \times \mc{S}_m$, which satisfies the Nehari and Poho\v{z}aev identities, to
\begin{equation*}
(-\Delta)^s u + \mu u = \overline g(u).
\end{equation*}
If $t_1 = +\infty$, then $\overline{g}(u_0) = g(u_0)$, so $(\mu_0,u_0)$ is a solution to \eqref{e-main}; otherwise, from Lemma \ref{lem_apriori_bound}, $|u_0|_\infty \le t_1$, hence $\overline{g}(u_0) = g(u_0)$ once again, and we conclude as before. In both cases, $\overline{G}(u_0) = G(u_0)$, where $\overline{G}' = \overline{g}$ and $\overline{G}(0) = 0$, so $K(u_0) < 0$.\\
Finally, if $g$ satisfies \eqref{eq_cond_meno_inf}, let $\omega>0$ to be determined and define $h_\omega$ and $H_\omega$ as in \eqref{e-h_omega}. 
Since $h_\omega(t_0) > 0$, let $t_{1,\omega}$ be the first zero of $h_\omega$ larger than $t_0$ ($t_{1,\omega} = +\infty$ if $h_\omega(t) > 0$ for all $t > t_0$) and set $\overline{h}_\omega(t) := h_\omega(t) \chi_{(0,t_{1,\omega})}(t)$.
In what follows, we will consider functionals and quantities analogous to those defined in \eqref{eq_def_Jmu}, \eqref{eq_def_Pmu}, and \eqref{eq_def_amu} replacing $g$ with $h_\omega$ or $\overline{h}_\omega$, using a subscript to distinguish the cases. 
Fix $\nu \in (0, \frac{G(t_0)}{t_0^2/2}) \subset (0,\overline{\mu}_{0})$. Arguing as in \cite[Proof of Theorem 2]{BeLi83I}, there exists $\phi \in H^s_\textup{rad}(\R^N)$ such that $0 \le \phi \le t_0$ and $\intRN G(\phi) - \frac{\nu}{2} \phi^2 \, \dx > 0$; then, in view of Remark \ref{r-all_mu}, up to internal scaling, we can assume that $\phi \in \mc{P}_{\nu}$. In particular, $H_{\omega}(\phi) = \overline{H}_{\omega}(\phi)$. Moreover, defining 
\begin{equation*}
\overline{\mu}_{\overline{h}_\omega} := \sup_{t \neq 0} \frac{\overline{H}_\omega(t)}{t^2/2},
\end{equation*}
we have 
$$\overline{\mu}_{\overline{h}_\omega} \geq \frac{\overline{H}_\omega(t_0)}{t_0^2/2} =\frac{G(t_0)}{t_0^2/2} + \omega > \nu + \omega.$$
In addition, $P_{\nu+\omega, \bar{h}_{\omega}}(\phi) = P_{\nu+\omega,h_{\omega}}(\phi) = P_{\nu}(\phi)=0$, thus
\begin{equation*}
\frac{1}{2}m_{\overline{h}_{\omega}} := \inf_{\nu' \in (0, \overline{\mu}_{\overline{h}_{\omega}})} \frac{a_{\omega,\overline{h}_{\omega}}(\nu')}{\nu'} \leq \frac{a_{\omega,\overline{h}_{\omega}}(\nu+\omega)}{\nu+\omega} \leq \frac{J_{\nu+\omega,\bar{h}_{\omega}}(\phi)}{\nu+\omega} = \frac{J_{\nu+\omega,h_{\omega}}(\phi)}{\nu+\omega} = \frac{J_{\nu}(\phi)}{\nu+\omega} \to 0
\end{equation*}
as $\omega \to +\infty$. We can then take $\omega$ so large that $m_{\overline{h}_\omega} < m$ and,
following the proof of Theorem \ref{thm_small_mass}, 
find a solution $(\mu_\omega,u_0) \in (0,+\infty) \times \mc{S}_m$ to
\begin{equation}\label{e-bar_h}
(-\Delta)^s u_0 + \mu_\omega u_0 = \overline{h}_\omega(u_0)
\end{equation}
that satisfies the Nehari and Poho\v{z}aev identities related to \eqref{e-bar_h}. Arguing as before, we obtain that, in fact, $(\mu_\omega,u_0)$ is a solution to
\begin{equation*}
(-\Delta)^s u_0 + \mu_\omega u_0 = h_\omega(u_0),
\end{equation*}
that is, $(\mu_\omega - \omega,u_0)$ is a solution to \eqref{e-main}. We conclude by letting $\mu_0 := \mu_\omega - \omega$ and observing that $(\mu_0,u_0)$ satisfies \eqref{eq_pohoz_frac} and \eqref{eq_nehar_id}.
\QED

%
%
%

\bigskip

We end this section with the 
examples mentioned in the introduction.

\medskip

\claim Proof of Corollary \ref{cor_log_power}.
As in \cite[Proof of Theorem 1.6]{MeSc24}, we prove that when $\beta < 0$ and $q \in (1,2^*_s-1]$, $\max G \lesseqgtr 0$ if and only if $\beta \lesseqgtr -\frac{\alpha (q+1)}{q-1} e^{-(q+1)/2}$, hence the second part (existence) follows from Theorem \ref{t-small-mass} 
and Proposition \ref{p-sym_unpert} (iv) (see also Remark \ref{r-ssm}). 
If $m>m_0$, then the claim 
follows by Theorem \ref{t-main}.
The first part (nonexistence) is a direct consequence of 
the Poho\v{z}aev identity, together with Proposition \ref{prop_intr_L2min_PN}.
%
%
\QED

\bigskip


\claim Proof of Corollary \ref{corol_low_power}.
We argue as in Corollary \ref{cor_log_power}.
\QED

\bigskip

\medskip

\textbf{Acknowledgments.}
This work was supported by the Thematic Research Programme ``Variational and geometrical methods in partial differential equations'', University of Warsaw, Excellence Initiative - Research University. The first author is supported by INdAM-GNAMPA Projects ``Metodi variazionali per problemi dipendenti da operatori frazionari isotropi e anisotropi'', codice CUP \#E5324001950001\#,
and ``Aspetti geometrici e qualitativi di equazioni ellitiche e paraboliche'', codice CUP \#E53C25002010001\#. The second author is supported by INdAM-GNAMPA Projects ``Problemi di ottimizzazione in PDEs da modelli biologici'', codice CUP \#E5324001950001\#, and ``Strutture Analitiche e Geometriche in PDEs: Regolarità, Fenomeni Critici e Dinamiche Complesse'', codice CUP \#E53C25002010001\#. The first author expresses his gratitude to the University of Warsaw for hosting him in November 2024 and February 2025.

\end{document}